\documentclass[leqno, 10pt]{amsart}
\usepackage{amsmath,amssymb,latexsym}
 \newtheorem{definition}{Definition}[section]
 \newtheorem{theorem}[definition]{Theorem}
 
 \newtheorem{proposition}[definition]{Proposition}
 \newtheorem{corollary}[definition]{Corollary}

 \newtheorem*{theorem*}{Theorem}
\newtheorem*{proposition*}{Proposition}
\newtheorem*{lemma*}{Lemma}

 \theoremstyle{remark}
 
 \newtheorem{remark}[definition]{Remark}

  \newtheorem*{acknowledgements}{Acknowledgements}
  \newtheorem*{notation}{Notation}

\newcommand{\op}[1]{\operatorname{#1}}

\newcommand{\tr}{\ensuremath{\op{tr}}}
\newcommand{\Tr}{\ensuremath{\op{Tr}}}
\newcommand{\Tra}{\ensuremath{\op{Trace}}}

\newcommand{\TR}{\ensuremath{\op{TR}}}

\newcommand{\Res}{\ensuremath{\op{Res}}}
\newcommand{\res}{\ensuremath{\op{res}}}

\newcommand{\ran}{\op{ran}}

\newcommand{\Sp}{\op{Sp}}

\newcommand{\End}{\ensuremath{\op{End}}}
\newcommand{\Cl}{\ensuremath{{\op{Cl}}}}

\newcommand{\C}{\ensuremath{\mathbb{C}}} 
 
\newcommand{\R}{\ensuremath{\mathbb{R}}} 
\newcommand{\Z}{\ensuremath{\mathbb{Z}}} 
\newcommand{\CZ}{\ensuremath{\mathbb{C}\!\setminus\!\mathbb{Z}}}

\newcommand{\Ca}[1]{\ensuremath{\mathcal{#1}}}

\newcommand{\cE}{\Ca{E}}

\newcommand{\cL}{\ensuremath{\mathcal{L}}}

\newcommand{\psido}{$\Psi$DO} 
\newcommand{\psidos}{$\Psi$DO's} 
 
\newcommand{\pdoi}{\ensuremath{\Psi^{\op{int}}}}

\newcommand{\ord}{{\op{ord}}}

\newcommand{\sD}{\ensuremath{{/\!\!\!\!D}}}

\newcommand{\sS}{\ensuremath{{/\!\!\!\!\!\;S}}}

\newcommand{\zetaup}{\zeta_{\scriptscriptstyle{\uparrow}}}
\newcommand{\zetadown}{\zeta_{\scriptscriptstyle{\downarrow}}}

\begin{document}
 
    \title{SPECTRAL ASYMMETRY, ZETA FUNCTIONS,\\ 
    AND THE NONCOMMUTATIVE RESIDUE}
\author{Rapha\"el Ponge}

\address{Department of Mathematics, Ohio State University, Columbus, USA.}
\curraddr{Max Planck Institute for Mathematics, Bonn, Germany.}
\email{raphaelp@mpim-bonn.mpg.de}
\thanks{Research partially supported by NSF grant DMS 0409005}

\numberwithin{equation}{section}

\begin{abstract}
In this paper we study the spectral asymmetry of (possibly nonselfadjoint) elliptic \psidos\ in terms of the difference of zeta functions coming 
from different cuttings. Refining previous formulas of Wodzicki in the case of odd class elliptic \psidos, our main results have several consequence concerning 
the local independence with respect to the cutting, the regularity at integer points of eta functions and a geometric expression for the spectral 
asymmetry of Dirac operators which, in particular, yields a new spectral interpretation of the Einstein-Hilbert action in gravity. 
\end{abstract}

\maketitle 

\section{Introduction}

This paper focuses on the spectral asymmetry of elliptic \psidos.  Given a compact Riemannian manifold $M^{n}$ and a Hermitian bundle $\cE$ over $M$, 
the spectral asymmetry was first  studied by 
 Atiyah-Patodi-Singer~\cite{APS:SARG1} in the case of a 
selfadjoint elliptic \psido\  $P:C^{\infty}(M,\cE) \rightarrow C^{\infty}(M,\cE)$ in terms of the eta function, 
\begin{equation}
    \eta(P;s)=\Tra P|P|^{-(s+1)}, \qquad s\in \C. 
\end{equation}
This function is meromorphic with at worst simple pole singularities and an important result, due to Atiyah-Patodi-Singer~\cite{APS:SARG3} and 
 Gilkey~(\cite{Gi:RGEFO}, \cite{Gi:ITHEASIT}), is its regularity at $s=0$,  so that the eta invariant $\eta(P):=\eta(P,0)$  
is always well defined. 

The residues of the eta function at other integer points are also interesting, e.g., they enter in the index formula of 
Br\"uning-Seeley~\cite{BS:ITFORSO} for first order elliptic operators on a manifold with cone-like singularities

In~\cite{Wo:SAZF}--\cite{Wo:PhD} Wodzicki took a different point of view. Motivated by an observation of Shubin,
he looked at the spectral asymmetry of  a (possibly nonselfadjoint) elliptic \psido\ 
$P:C^{\infty}(M,\cE) \rightarrow C^{\infty}(M,\cE)$ of order $m>0$ in terms of the difference, 
 \begin{equation}
    \zeta_{\theta}(P;s)-\zeta_{\theta'}(P;s)=\Tra P_{\theta}^{-s}-\Tra P_{\theta'}^{-s}, \qquad s \in \C,
 \end{equation}
    of zeta functions  coming from different spectral cuttings $L_{\theta}=\{\arg\lambda =\theta\}$ and $L_{\theta'}=\{\arg\lambda =\theta'\}$ with 
 $0\leq \theta<\theta'<2\pi$.  In particular, he showed that the spectral asymmetry of $P$ was encoded by the sectorial projection earlier introduced by 
 Burak~\cite{Bu:SPEO} and given by 
   \begin{equation}
       \Pi_{\theta,\theta'}(P)=\frac{1}{2i\pi} \int_{\Gamma_{\theta,\theta'}}\lambda^{-1}P(P-\lambda)^{-1}d\lambda, 
       \label{eq:Intro.sectorial-projection}
     \end{equation}
where $\Gamma_{\theta,\theta'}$ is a contour separating the part of the spectrum of $P$ contained in the open sector 
$\theta<\arg \lambda< \theta'$ from the rest of the spectrum. More precisely, Wodzicki
proved the equality of meromorphic functions, 
    \begin{equation}
     \zeta_{\theta}(P;s)- \zeta_{\theta'}(P;s)=(1-e^{-2i\pi s}) \Tra \Pi_{\theta,\theta'}(P)P_{\theta}^{-s}, \quad s 
     \in \C. 
     \label{eq:intro-Wodzicki-Formula1}
   \end{equation}
In particular, at every integer $k\in \Z$ the function $ \zeta_{\theta}(P;s)- \zeta_{\theta'}(P;s)$ is regular  and there we 
have 
   \begin{equation}
     \ord P. \lim_{s\rightarrow k}(\zeta_{\theta}(P;s)-\zeta_{\theta'}(P;s))=2i\pi \Res 
\Pi_{\theta,\theta'}(P)P^{-k}, 
       \label{eq:intro-Wodzicki-Formula2}
   \end{equation}
where $\Res$ denotes the noncommutative residue of Wodzicki~(\cite{Wo:LISA}, \cite{Wo:NCRF}) and 
Guillemin~\cite{Gu:NPWF}.

Furthermore, Wodzicki proved in~\cite{Wo:LISA} that the regular value $\zeta(P;0)$ is independent of the choice of the cutting $L_{\theta}$ 
and that  the noncommutative residue of a \psido\ projection is always zero,  which generalize the vanishing of the residue at the 
origin of the eta function a selfadjoint elliptic \psido. 
%


In this paper, partly motivated by a recent upsurge of interest in the spectral asymmetry of non-selfadjoint elliptic \psidos~(\cite{BK:RAT}, 
\cite{Sc:RD}), 
we prove various results related to  
the spectral asymmetry of odd class elliptic \psidos\ as a consequence of a refinement of 
the formulas~(\ref{eq:intro-Wodzicki-Formula1})--(\ref{eq:intro-Wodzicki-Formula2}) for such operators.

Recall that a \psido\ of integer order is said to be odd class when the homogeneous 
components  of its symbol are homogeneous with respect to the dilation by $-1$. In particular,  the odd class \psidos\ form an algebra containing all the differential 
 operators and the parametrices of odd class elliptic \psidos. 
 
 Let $P:C^{\infty}(M,\cE) \rightarrow C^{\infty}(M,\cE)$ be an odd class \psido\ of integer order $m\geq 1$ and 
 let $L_{\theta}=\{\arg\lambda =\theta\}$ and $L_{\theta'}=\{\arg\lambda =\theta'\}$ be spectral cuttings for $P$ and its principal symbol with 
 $0\leq \theta<\theta'<2\pi$. Our main results are:\smallskip 
 
   (i) If $\dim M$ is odd and $\ord P$ is even then $\zeta_{\theta}(P;s)$ is regular at every integer point and its 
        value there is independent of the spectral cut $L_{\theta}$ (Theorem~\ref{SAZF.odd-PsiDOs1}).\smallskip
        
  (ii)  If $\dim M$ is even, $\ord P$ is odd and the principal symbol of $P$ has all its eigenvalue in the open 
      cone $\{\theta <\arg \lambda<\theta'\}\cup\{\theta+\pi <\arg \lambda<\theta'+\pi\}$, then for any integer $k\in \Z$ we have  
      \begin{equation}
          \ord P.\lim_{s\rightarrow k}(\zeta_{\theta}(P;s)- \zeta_{\theta'}(P;s))=i\pi \Res P^{-k}.
      \end{equation}
      In particular, at every integer at which they are not singular the  functions 
      $\zeta_{\theta}(P;s)$ and $\zeta_{\theta'}(P;s)$ take on the same regular value (Theorem~\ref{SAZF.odd-PsiDOs2}).\smallskip
      
 These results are deduced from a careful analysis of the symbol of the sectorial projection $\Pi_{\theta,\theta'}(P)$, so that the proofs are 
 purely local in nature. It thus follows that the theorems ultimately hold at the level of the \emph{local} 
 zeta functions $\zeta_{\theta}(P;0)(x)$ and $\zeta_{\theta'}(P;0)(x)$,  
 that is, the densities whose integrals yield the zeta functions $\zeta_{\theta}(P;s)$ and $\zeta_{\theta'}(P;s)$. In particular, we obtain that if $P$ is an odd class elliptic 
 \psido\ 
 sastisfying either the assumptions of (i) or that of (ii), then the regular value $\zeta_{\theta}(P;0)(x)$ is independent of the choice of the 
spectral cutting (Theorem~\ref{SAZF.odd-PsiDOs-local}). 

In fact, the independence with respect to the spectral cutting of the regular values at $s=0$ of the local zeta functions is 
not true for general \psidos\ (see~\cite[pp.~130-131]{Wo:SAZF}). Therefore, it is interesting to see that this nevertheless can happen for a wide class of elliptic 
\psidos. 

Next, these results have further applications when $P$ is selfadjoint. In this case we shall 
 use the subscript $\uparrow$  (resp.~$\downarrow$) to refer to a spectral cut in the upper halfplane $\Im \lambda>0$ 
(resp.~lower halfplane $\Im \lambda<0$).

First, while the above results tell us that there are many integer points at which 
there is no spectral asymmetry, they also allow us to single out some points at which the spectral asymmetry always occurs. For instance, we always have 
\begin{equation}
    \lim_{s\rightarrow n}\frac{1}{i}(\zetaup(P;s)- \zetadown(P;s))> 0, 
     \label{eq:Intro.asymmetry-point}
\end{equation}
when $\dim M$ is even and $P$ is a first order selfadjoint elliptic odd class \psido\ (see Proposition~\ref{prop:SAZF.selfadjoint-asymmetry}).

Second, as the eta function $\eta(P;s)$ can be nicely related to $\zetaup(P;s)-\zetadown(P;s)$ (see \cite[p.~114]{Sh:POST} and 
Section~\ref{sec.selfadjoint-odd}),   
we can make use of the previous results to study $\eta(P;s)$. It is a well known result of Branson-Gilkey~\cite{BG:REFODT} that in even 
dimension the eta function of a Dirac operator is an entire function. We generalize this result by proving that 
if $\ord P$ and $\dim M$ have opposite parities then $\eta(P;s)$ is regular at every integer point, so that when 
$P$ has order 1 and  $\dim M$ is even the function $\eta(P;s)$ is entire (Theorem~\ref{thm:SAZF.regularity-eta}).

The latter  result has been independently obtained by Grubb~\cite{Gr:RATZLE} using a different approach.  Furthermore, it allows us
to simplify in odd dimension the aforementioned index formula of Br\"uning-Seeley~\cite{BS:ITFORSO} 
for first order elliptic operators on a manifold with cone-like singularities (see Remark~\ref{rem:SAZF.Bruening-Seeley}).

Third, for Dirac operators our results enable us to express the spectral asymmetry of these operators in geometric 
terms. More precisely, assume that $M$ has even dimension, that $\cE$ is a  $\Z_{2}$-graded Clifford module 
over $M$ equipped with a unitary connection $\nabla^{\cE}$, and let $\sD_{\cE}: C^{\infty}(M,\cE) \rightarrow C^{\infty}(M,\cE)$ be the associated Dirac 
operator.    
Then in Proposition~\ref{thm:SAZF.Dirac-operator} we show that:\smallskip

- At every  integer that is not an even integer between $2$ and $n$ the zeta functions 
     $\zetaup(\sD_{\cE};s)$ and  $\zetadown(\sD_{\cE};s)$ are non-singular and take on the same regular value;
     
-   For  $k=2,4,\ldots,n$ we can express $  \lim_{s\rightarrow k}(\zetaup(\sD_{\cE};s)- \zetadown(\sD_{\cE};s))$ as the integral of 
     a universal polynomial in complete tensorial contractions of the covariant derivatives of the  
     curvature $R^{M}$ of $M$ and the twisted curvature $F^{\cE/\sS}$ of $\cE$.
 
As a consequence we  get a new spectral interpretation of the Einstein-Hilbert action $\mathcal{I}=
\int_{M}r_{M}(x)\sqrt{g(x)}dx$, which is an important issue in noncommutative geometry and we get points at which the spectral asymmetry occurs 
indepently of the choice of the Clifford data $(\cE, \nabla^{\cE})$  (see Proposition~\ref{prop:SAZF.Dirac-asymmetry} ). 

The paper is organized as follows. In Section~\ref{sec.CPNCR} we recall the general background needed in this paper about complex powers of elliptic operators, the 
noncommutative residue trace of Wodzicki and Guillemin and the zeta and eta functions of elliptic \psidos. In Section~\ref{sec.SP} we gather some of the main facts 
about the sectorial projection of an elliptic \psido, but we postpone to the Appendix those concerning its spectral interpretation. 
In Section~\ref{sec:SAF} we give a detailed review of Wodzicki's results on 
the spectral asymmetry elliptic \psidos\ needed in this paper. 
In Section~\ref{sec.odd} we refine the latter formulas for odd class elliptic \psidos\ and prove our main results. 
We then specialize these results to the selfadjoint case in Section~\ref{sec.selfadjoint-odd} and to Dirac operators in Section~\ref{sec.Dirac}. 

\begin{notation}
 Throughout all this paper we let $M$ denote a compact Riemannian manifold of dimension $n$ and  let $\cE$ 
be a Hermitian vector bundle over $M$ of rank $r$. 
  \end{notation}

\section{General background}
\label{sec.CPNCR}

In this section we recall the main facts about
complex powers of elliptic \psidos,  the  
noncommutative residue trace of Wodzicki and Guillemin and the  zeta and eta functions of elliptic \psidos.
\subsection{Complex powers of elliptic \psidos}
For $m\in \C$ we let $\Psi^{m}(M,\cE)$ denote the space of (classical) \psidos\ of order $m$ on $M$ acting on sections of $\cE$, 
i.e.,~continuous operators $P:C^{\infty}(M,\cE) \rightarrow C^{\infty}(M,\cE)$ such that:\smallskip 

- The distribution kernel of $P$ is smooth off the diagonal of 
$M\times M$;\smallskip 

- In any local trivializing chart  $U \subset \R^{n}$ the operator $P$ is of the form 
$P=p(x,D)+R$, for some polyhomogeneous symbol $p(x,\xi)\sim \sum_{j \geq 0}p_{m-j}(x,\xi)$  of degree $m$ and some smoothing operator $R$, where  
$p(x,D)$ denotes the linear operator from $C_{c}^{\infty}(U, \C^{r})$ to $C^{\infty}(U, \C^{r})$ such that 
\begin{equation}
    p(x,D)u(x)=(2\pi)^{-n}\int e^{ix.\xi}p(x,\xi) \hat{u}(\xi)d\xi \quad \forall u \in C^{\infty}_{c}(U, \C^{r}).
\end{equation}

Let $P:C^{\infty}(M,\cE)\rightarrow C^{\infty}(M,\cE)$ be an elliptic \psido\ of degree $m>0$ with principal 
symbol $p_{m}(x,\xi)$ and assume that the ray $L_{\theta}=\{\arg \lambda=\theta\}$, $0\leq \theta<2\pi$, is a spectral cutting for $p_{m}$, 
that is, $p_{m}(x,\xi)-\lambda$ 
is invertible for every $\lambda \in L_{\theta}$. Then there is a conical neighborhood  $\Lambda$ of $L_{\theta}$ such that any ray contained in 
$\Lambda$ is also a spectral cutting for $p_{m}$. It then follows that $P$ admits an asymptotic resolvent as a parametrix in a suitable class of \psidos\ with 
parametrized by $\Lambda$ (see \cite{Se:CPEO}, \cite{Sh:POST}, \cite{GS:WPPOAPSBP}). This allows us to show that, for any closed cone 
$\Lambda'$ such that $\overline{\Lambda'}\setminus 0\subset \Lambda$ and for $R>0$ large enough, there 
exists $C_{\Lambda' R}>0$ such that 
\begin{equation}
    \|(P-\lambda)^{-1}\|_{\cL(L^{2}(M,\cE))}\leq C_{\Lambda' R}|\lambda|^{-1}, \qquad \lambda \in \Lambda', \quad |\lambda|\geq R.
     \label{eq:NCR.minimal-growth}
\end{equation}
Therefore there are infinitely many rays $L_{\theta}=\{\arg \lambda =\theta\}$ contained in $\Lambda$ that are not through an eigenvalue of 
$P$ and any such ray is a ray of minimal growth. 
 
 On the other hand, (\ref{eq:NCR.minimal-growth}) also implies that the spectrum of $P$ is not $\C$, hence
 consists of an unbounded set of isolated eigenvalues with finite multiplicities. Thus, 
 we can define the root space and Riesz projection associated to $\lambda \in \Sp P$  by letting 
\begin{equation}
    E_{\lambda}(P)=\cup_{j \geq 1}\ker(P-\lambda)^{j} \quad \text{and} \quad \Pi_{\lambda}(P)=\frac{-1}{2i\pi} 
     \int_{\Gamma_{(\lambda)}} (P-\mu)^{-1}d\mu,
\end{equation}
where $\Gamma_{(\lambda)}$ is a direct-oriented circle about $\lambda$ with a radius small enough so that 
apart from $\lambda$ no 
other element of $\Sp P \cup \{0\}$  lies inside $\Gamma_{(\lambda)}$. 
 
The family $\{\Pi_{\lambda}(P)\}_{\lambda \in \Sp P}$ is a family of disjoint projections, 
in the sense that we have $\Pi_{\lambda}(P)\Pi_{\mu}(P)=0$ for $\lambda\neq \mu$. Moreover, for every $\lambda \in \Sp P$ the root space 
$E_{\lambda}(P)$ has finite dimension and 
$\Pi_{\lambda}(P)$ projects onto $E_{\lambda}(P)$ 
and along $E_{\bar{\lambda}}(P^{*})^{\perp}$ (see~\cite[Le\c con 148]{RN:LAF}, \cite[Sect.~I.7]{GK:ITLNSO}). In addition, since $P$ is elliptic 
$\Pi_{\lambda}(P)$ is a smoothing operator and $E_{\lambda}(P)$ is contained  
$C^{\infty}(M,\cE)$ (see \cite[Thm.~8.4]{Sh:POST}). 
 
 Next, assume that the ray $L_{\theta}$ is a spectral cutting for both $p_{m}$ and $P$. 
 Then the family $(P_{\theta}^{s})_{s \in \C}$ of complex powers of $P$ associated to $L_{\theta}$ can be defined as follows. Thanks 
 to~(\ref{eq:NCR.minimal-growth}) we define a bounded operator on $L^{2}(M,\cE)$ by letting
\begin{gather}
P_{\theta}^{s}= \frac{-1}{2i\pi} \int_{\Gamma_{\theta}} \lambda^{s}_{\theta}(P-\lambda)^{-1}d\lambda, \qquad \Re s<0,
     \label{eq:background.complex-powers-definition}\\
\Gamma_{\theta}=\{ \rho e^{i\theta}; \infty <\rho\leq r\}\cup\{ r e^{it}; 
\theta\geq t\geq \theta-2\pi \}\cup\{ \rho e^{i(\theta-2\pi)};  r\leq \rho\leq \infty\},  
\label{eq:background.complex-powers-definition-Gammat}
\end{gather}
where $r>0$ is small enough so that there is no nonzero eigenvalue of $P$ in the disc $|\lambda|<r$ and 
$\lambda^{s}_{\theta}=|\lambda|^{s}e^{is\arg_{\theta}\lambda}$ is defined by means of the continuous determination of the argument on $\C\setminus 
L_{\theta}$ that takes values in $(\theta-2\pi,\theta)$. We then have
\begin{gather}
    P_{\theta}^{s_{1}+s_{2}}=P_{\theta}^{s_{1}} P_{\theta}^{s_{2}}, \quad \Re s_{j}<0,\\
    P_{\theta}^{-k}=P^{-k}, \quad k=1,2,\ldots, 
\end{gather}
where $P^{-k}$ denote the partial inverse of $P^{k}$, that is, the bounded operator that inverts $P$ on 
$E_{0}(P^{k*})^{\perp}=E_{0}(P^{*})^{\perp}$ and vanishes on $E_{0}(P^{k})=E_{0}(P)$. 

On the other hand, the \psido\ calculus with parameter allows us to show 
that $P_{\theta}^{s}$ is a \psido\ of order $ms$ and that the family $(P_{\theta}^{s})_{\Re s<0}$ is  a holomorphic family of \psidos\ in the 
sense of~\cite[7.14]{Wo:LISA} and~\cite[p.~189]{Gu:GLD} (see~\cite{Se:CPEO}, \cite{Sh:POST}, \cite{GS:WPPOAPSBP}). 
Therefore, for any $s \in \C$ we can define $P_{\theta}^{s}$ as the \psido\ such that $P^{s}_{\theta}=P^{k} P^{s-k}_{\theta}$, where $k$ is any integer~$>\Re s$.  

This gives rise to a holomorphic 1-parameter group of \psidos\ such that $\ord P_{\theta}^{s}=ms$ for any $s\in \C$. In particular, we have 
$P_{\theta}^{0}=PP^{-1}=1-\Pi_{0}(P)$. 

\subsection{Noncommutative residue} 
The noncommutative residue trace of Wodzicki~(\cite{Wo:LISA}, \cite{Wo:NCRF}) and Guillemin~\cite{Gu:NPWF} appears as the residual trace 
on the algebra $\Psi^{\Z}(M,\cE)$  of \psidos\ of integer orders induced by the analytic extension of the usual trace to the class 
$\Psi^{\CZ}(M,\cE)$  of \psidos\ of non-integer complex orders. Our exposition essentially follows that of~\cite{KV:GDEO} 
and~\cite{CM:LIFNCG}.

First, if $Q$ is in $\pdoi(M,\cE) = \cup_{\Re m < -n}\Psi^{m}(M,\cE)$ then the restriction of its distribution kernel  
to the diagonal of $M\times M$  an element $k_{Q}(x,x)$ of $\Gamma(M,|\Lambda|(M)\otimes \End \cE)$, the space of smooth  $\End \cE$-valued 
densities. Therefore, the operator $Q$ is trace-class and we have 
$ \Tra Q = \int_{M} \tr_{\cE}k_{Q}(x,x)$.

In fact, as shown in~\cite{KV:GDEO} the map $Q \rightarrow k_{Q}(x,x)$ has a unique analytic 
continuation $Q \rightarrow t_{Q}(x)$ to the class $\Psi^{\CZ}(M,\cE)$, where analyticity is meant in  
in the sense that, for every 
holomorphic family $(Q_{z})_{z \in \Omega}$ with values in $\Psi^{\CZ}(M,\cE)$, the map $z \rightarrow t_{Q_{z}}(x)$ is analytic with 
values in $\Gamma(M,|\Lambda|(M)\otimes \End \cE)$.

Moreover, if $Q$ is in $\Psi^{\Z}(M,\cE)$ and $(Q_{z})_{z\in \Omega}$ is a holomorphic family of \psidos\ defined near $z=0$ such that $Q_{0}=Q$ and 
$\ord Q_{z}=z+\ord Q$, then the map $z\rightarrow t_{Q_{z}}(x)$ has at worst a simple pole singularity at $z=0$  
in such way that in local trivializing coordinates we have 
\begin{equation}
\res_{z=0}t_{Q_{z}}(x)= - (2\pi)^{-n}\int_{|\xi|=1}q_{-n}(x,\xi) d^{n-1}\xi,
\end{equation}
where $q_{-n}(x,\xi)$ denotes the symbol of degree $-n$ of $Q$. Since $t_{Q_{z}}(x)$ is a density we see that we 
get a well defined $\End \cE$-valued density on $M$ by letting
\begin{equation}
    c_{Q}(x)=  (2\pi)^{-n}(\int_{|\xi|=1}q_{-n}(x,\xi) d^{n-1}\xi) .
     \label{eq:NCR.density}
\end{equation}

We can now define the functionals 
\begin{gather}
     \TR Q= \int_{M} \tr_{\cE} t_{Q}(x), \qquad Q\in \Psi^{\CZ}(M,\cE),\\
 \Res Q = \int_{M} \tr_{\cE} c_{Q}(x), \qquad Q\in \Psi^{\Z}(M,\cE).
\end{gather}

\begin{theorem}[\cite{KV:GDEO}]\label{thm:background.TR}
   1) The functional $\TR$  is the unique analytic continuation of the usual trace to $\Psi^{\CZ}(M,\cE)$.\smallskip 
   
   2) We have $\TR [Q_{1},Q_{2}]=0$ whenever $\ord Q_{1}+\ord Q_{2}\not\in \Z$. \smallskip 
   
   3) Let $Q \in \Psi^{\Z}(M,\cE)$ and let $(Q_{z})_{z\in \Omega}$ be a holomorphic family of \psidos\ 
defined near $z=0$ such that $Q_{0}=Q$ and $\ord Q_{z}=z+\ord Q$. Then near $z=0$ 
the function $\TR Q_{z}$ has at worst a simple pole singularity such  that $\res_{z=0}\TR Q_{z}=-\Res Q$.
\end{theorem}

The functional $\Res$ is the noncommutative residue of Wodzicki and Guillemin. From Theorem~\ref{thm:background.TR} we immediately get:

\begin{theorem}[\cite{Wo:LISA}, \cite{Gu:NPWF}, \cite{Wo:NCRF}]\label{prop:background.NCR}
   1) The noncommutative residue is a linear trace on the algebra $\Psi^{\Z}(M,\cE)$ which 
   vanishes on differential operators and on \psidos\ of integer order $\leq -(n+1)$.\smallskip 
   
   2) We have $\res_{s=0} \TR QP_{\theta}^{-s}=m \Res Q$ for any $Q\in \Psi^{\Z}(M,\cE)$.
\end{theorem}

Notice also that by a well-known result of Wodzicki~(\cite{Wo:PhD}, \cite[Prop.~5.4]{Ka:RNC};
see also~\cite{Gu:RTCAFIO}) if $M$ is connected and 
has dimension $\geq 2$ then the noncommutative residue induces the only trace on
$\Psi^{\Z}(M,\cE)$ up to a multiplicative constant. 

\subsection{Zeta and eta functions} 
The canonical trace $\TR$ allows us to define the zeta function of $P$ 
as the meromorphic function on $\C$ given  by 
\begin{equation}
    \zeta_{\theta}(P;s)=\TR P_{\theta}^{-s}, \quad s\in \C. 
\end{equation}
Then from Theorem~\ref{thm:background.TR}  we obtain: 

\begin{proposition}\label{prop:background.zeta-function} Let $\Sigma=\{\frac{n-j}{m};\ j=0,1, \ldots\}\setminus \{0\}$. Then  $\zeta_{\theta}(P;s)$ 
    is analytic outside $\Sigma$ and on $\Sigma$ has at  worst simple pole singularities such that 
 \begin{equation}
     \res_{s=\sigma}\zeta_{\theta}(P;s)=m\Res P^{-\sigma}_{\theta}, \qquad \sigma\in \Sigma.
     \label{eq:Background.residues-zeta-NCR}
 \end{equation}
\end{proposition}

Notice that~(\ref{eq:Background.residues-zeta-NCR}) is true for $\sigma=0$ as well, but in this case it gives 
\begin{equation}
    \res_{s=0}\zeta_{\theta}(P;s)=\Res P^{0}_{\theta}=\Res 
    [1-\Pi_{0}(P)]=0,
\end{equation}
 since $\Pi_{0}(P)$ is a smoothing operator. Thus $\zeta_{\theta}(P;s)$ is always regular at $s=0$. 

Finally, assume that $P$ is selfadjoint. Then  the eta function of $P$ is the meromorphic function given by
\begin{equation}
    \eta(P;s)=\TR F|P|^{-s}, \qquad s \in \C, 
\end{equation} 
where $F=P|P|^{-1}$ is the sign operator of $P$. Then  using
Theorem~\ref{thm:background.TR}  we get:  

\begin{proposition}\label{prop:background.eta-function}Let $\Sigma=\{\frac{n-j}{m};\ j=0,1, \ldots\}$. Then 
$\eta(P;s)$ is analytic outside $\Sigma$ and on $\Sigma$ has at  worst simple pole singularities such that 
 \begin{equation}
     \res_{s=\sigma}\eta(P;s)=m\Res  F|P|^{-\sigma},  \qquad \sigma\in \Sigma.
     \label{eq:background.residues-eta-function}
 \end{equation}
\end{proposition}

Showing the regularity at the origin of $\eta(P;s)$ is a much more difficult task than for the zeta 
functions. Indeed, from~(\ref{eq:background.residues-eta-function}) we get 
\begin{equation}
     \res_{s=0}\eta(P;s)=m\Res F=m\int_{M}\tr_{\cE}c_{F}(x),
\end{equation}
 and examples 
show that  $c_{F}(x)$ need not vanish locally (see \cite{Gi:RLEFO}). Therefore, Atiyah-Patodi-Singer~\cite{APS:SARG3} and 
    Gilkey~(\cite{Gi:RGEFO}, \cite{Gi:ITHEASIT}) had to rely on global and $K$-theoretic arguments to prove: 

\begin{theorem}
    The function $\eta(P;s)$ is always regular at $s=0$. 
\end{theorem}
This shows that the eta invariant 
$\eta(P):= \eta(P;0)$ is always well defined. Since its appearance as a boundary correcting term 
in the index formula of Atiyah-Patodi-Singer~\cite{APS:SARG1}, the eta invariant has found many applications and has 
been extended to various other settings. We refer to the surveys of Bismut~\cite{Bi:LITEIHTS} and 
M\"uller~\cite{Mu:EI}, and the references therein, for an overview of the main results on the eta invariant. 

\section{The sectorial projection of an elliptic \psido}
\label{sec.SP}
In this section we give a detailed account on the sectorial projection of an elliptic \psido\ introduced by Burak~\cite{Bu:SPEO}. 
 
Let $P:C^{\infty}(M,\cE)\rightarrow 
C^{\infty}(M,\cE)$ be an elliptic \psido\ of order $m>0$ 
and assume that $L_{\theta}=\{\arg \lambda =\theta\}$  and 
$L_{\theta'}=\{\arg \lambda =\theta\}$ are spectral cuttings for both $P$ and its principal 
symbol $p_{m}(x,\xi)$ with $\theta< \theta'\leq \theta+2\pi$. In addition,
we let $\Lambda_{\theta,\theta'}$ and $\Lambda_{\theta',\theta+2\pi}$ respectively denote the angular sectors $\theta<\arg \lambda<\theta'$ 
and $\theta'<\arg \lambda<\theta+2\pi$. 

The sectorial projection of $P$ associated to the angular sector $\Lambda_{\theta,\theta'}$ is 
  \begin{gather}
       \Pi_{\theta,\theta'}(P)=\frac{1}{2i\pi} \int_{\Gamma_{\theta,\theta'}}\lambda^{-1}P(P-\lambda)^{-1}d\lambda, 
         \label{eq:SP.SAPi}\\
      \Gamma_{\theta,\theta}=\{ \rho e^{i\theta}; \infty >\rho\geq r\} \cup \{r e^{it}; \theta\leq t\leq \theta'\} 
 \cup  \{ \rho e^{i\theta'}; r\leq \rho< \infty\},   \label{eq:SP.SAPi.Gammatt'} 
     \end{gather}
 where $r$  is small enough so that no non-zero eigenvalue of $P$ lies in the disc $|\lambda|\leq r$.

In view of~(\ref{eq:NCR.minimal-growth}) the integral~(\ref{eq:SP.SAPi}) \emph{a priori} gives rise to an unbounded operator on $L^{2}(M,\cE)$ 
whose domain contains $L_{m}^{2}(M,\cE)$. We actually get a bounded operator thanks to: 

\begin{proposition}\label{prop:Appendix.sectorial-projection-PsiDO}
  1)  The operator $\Pi_{\theta,\theta'}(P)$ is a \psido\ of order~$\leq 0$,  hence is bounded on $L^{2}(M,\cE)$.\smallskip
    
   2) The zero'th order symbol of $\Pi_{\theta,\theta'}(P)$ is the sectorial projection $\Pi_{\theta,\theta'}(p_{m}(x,\xi))$, i.e.,~the Riesz 
   projection onto the root space associated to eigenvalues in $\Lambda_{\theta,\theta'}$. 
\end{proposition}
\begin{proof}
Let $R_{\theta,\theta'}=\frac{1}{2i\pi} \int_{\Gamma_{\theta,\theta'}}\lambda^{-1}(P-\lambda)^{-1}d\lambda$. 
Then the arguments of~\cite[Thm.~3]{Se:CPEO} can be carried through to prove that $R_{\theta,\theta'}$ 
is a \psido\ of order $\leq -1$. Hence $\Pi_{\theta,\theta'}(P)=P R_{\theta,\theta'}$ is a \psido\ of order~$\leq 0$.
 
 Next, in some local trivializing coordinates let $p(x,\xi)\sim \sum_{j\geq 0}p_{m-j}(x,\xi)$ and $r(x,\xi) \sim \sum_{j\geq 0} r_{-1-j}(x,\xi)$ 
 respectively denote the symbols of $P$ 
 and $R_{\theta,\theta'}$, so that $\Pi_{\theta,\theta'}(P)$ has symbol $\pi(x,\xi)\sim \sum 
 \frac{(-i)^{|\alpha|}}{\alpha!}\partial_{\xi}^{\alpha}p(x,\xi)\partial_{x}^{\alpha}r(x,\xi)$. Furthermore, let 
 $q(x,\xi) \sim \sum_{j\geq 0} q_{-m-j}(x,\xi;\lambda)$ be the symbol with parameter of $(P-\lambda)^{-1}$. 
 Then by~\cite[Thm.~2]{Se:CPEO} we have 
 \begin{equation}
     r(x,\xi) = \frac{1}{2i\pi} \int_{\Gamma_{\theta,\theta'}}\lambda^{-1}q(x,\xi;\lambda)d\lambda= 
     \frac{-1}{2i\pi} \int_{\Gamma_{(x,\xi)}}\lambda^{-1}q(x,\xi;\lambda)d\lambda,
 \end{equation}
 where $\Gamma_{(x,\xi)}$ is a direct-oriented bounded contour contained in the sector $\Lambda_{\theta,\theta'}$ 
 which isolates from $\C\setminus \Lambda_{\theta,\theta'}$ the eigenvalues of $p_{m}(x,\xi)$ that lie in 
 $\Lambda_{\theta,\theta'}$. 
 
 On the other hand, using the equality, 
 \begin{equation}
   P(P-\lambda)^{-1}=1+\lambda(P-\lambda)^{-1},
    \label{eq:Appendix.proof.trick}
\end{equation}
we see that $\lambda^{-1}(p\# q)(x,\xi;\lambda)=\lambda^{-1}+q(x,\xi;\lambda)$. Thus $\pi(x,\xi)$ is equal to
\begin{equation}
    p\#r(x,\xi)=  \frac{-1}{2i\pi} \int_{\Gamma_{(x,\xi)}}\lambda^{-1}p\#q(x,\xi;\lambda)d\lambda =    \frac{-1}{2i\pi} 
    \int_{\Gamma_{(x,\xi)}}q(x,\xi;\lambda)d\lambda.  
\end{equation}
Therefore, for $j=0,1,\ldots$ we obtain
 \begin{equation}
     \pi_{-j}(x,\xi)= \frac{-1}{2i\pi} \int_{\Gamma_{(x,\xi)}} q_{-m-j}(x,\xi;\lambda) d\lambda.
     \label{eq:Appendix.symbol-Wodzicki-projection}
\end{equation}
 Hence $\pi_{0}(x,\xi)= \frac{-1}{2i\pi} \int_{\Gamma_{(x,\xi)}} (p_{m}(x,\xi)-\lambda)^{-1} d\lambda=\Pi_{\theta,\theta'}(p_{m}(x,\xi))$ as 
 desired.\end{proof}

Next, the sectorial root spaces $E_{\theta,\theta'}(P)$ and $E_{\theta',\theta+2\pi}(P)$ are 
\begin{equation}
    E_{\theta,\theta'}(P)= \dotplus_{\lambda \in \Lambda_{\theta,\theta'}} E_{\lambda}(P),  \qquad 
    E_{\theta',\theta+2\pi}(P)= \dotplus_{\lambda \in \Lambda_{\theta',\theta+2\pi}} E_{\lambda}(P),
\end{equation}
where $\dotplus$ denotes the algebraic direct sum and for $\lambda \not \in \Sp P$ we make the convention that  
$E_{\lambda}(P)=\cup_{k\geq 1}\ker (P-\lambda)^{k}=\{0\}$. Then we have: 

\begin{proposition}\label{prop:SP.range-kernel}
    $\Pi_{\theta,\theta'}(P)$ is a projection on $L^{2}(M,\cE)$ which projects onto a subspace containing $E_{\theta,\theta'}(P)$ and 
    along a subspace containing $E_{0}(P)\dotplus 
E_{\theta',\theta+2\pi}(P)$.
\end{proposition}
\begin{proof}
Let $L_{\theta_{1}}$ and $L_{\theta_{2}}$ be rays with $\theta_{1}<\theta<\theta'<\theta_{1}'<\theta+2\pi$ and such that no eigenvalues of $P$ and $p_{m}$ 
lie in the angular sectors $\theta_{1}<\arg \lambda<\theta$ and $\theta'<\arg \lambda <\theta$. This allows us to replace in the formula~(\ref{eq:SP.SAPi}) for 
$\Pi_{\theta,\theta'}(P)$  the integration over $\Gamma_{\theta,\theta'}$ by that over a contour $\Gamma_{\theta_{1},\theta'_{1}}$ defined 
as in~(\ref{eq:SP.SAPi.Gammatt'}) using $\theta_{1}$ and $\theta_{1}'$ and a radius $r_{1}$ smaller than that of $\Gamma_{\theta,\theta'}$. Then we 
have
\begin{equation}
   \Pi_{\theta,\theta'}(P)^{2} = \frac{-1}{4\pi^{2}} \int_{\Gamma_{\theta,\theta'}}\int_{\Gamma_{\theta_{1},\theta_{1}'}} 
        \lambda^{-1}\mu^{-1}P^{2}(P-\lambda)^{-1}(P-\mu)^{-1}d\lambda d\mu.
\end{equation}
Therefore, by using the 
identity,
\begin{equation}
    (P-\lambda)^{-1}(P-\mu)^{-1}=(\lambda-\mu)^{-1}[(P-\lambda)^{-1}-(P-\mu)^{-1}],
     \label{eq:Appendix-spectral.identity1}
\end{equation}
we deduce that $-4\pi^{2}\Pi_{\theta,\theta'}(P)^{2}$ is equal to
\begin{equation}
  \label{eq:SAZF.Pitt'-projection}
     \int_{\Gamma_{\theta,\theta'}} \frac{P^{2}}{\lambda(P-\lambda)}
        (\int_{\Gamma_{\theta_{1},\theta_{1}'}}\frac{\mu^{-1}d\mu}{\mu-\lambda})d\lambda +       
       \int_{\Gamma_{\theta_{1},\theta'_{1}}}\frac{P^{2}}{\mu(P-\mu)} 
        (\int_{\Gamma_{\theta,\theta}}\frac{\lambda^{-1}d\mu}{\lambda-\mu})d\mu,
\end{equation}
from which we see that $\Pi_{\theta,\theta'}(P)^{2}= \frac{1}{2i\pi} \int_{\Gamma_{\theta,\theta'}} \lambda^{-2}P^{2}(P-\lambda)^{-1}d\lambda$. 
Combining this with~(\ref{eq:Appendix.proof.trick}) 
then gives
\begin{equation}
    \Pi_{\theta,\theta'}(P)^{2}=  \frac{1}{2i\pi} \int_{\Gamma_{\theta,\theta'}}\lambda^{-2}Pd\lambda +  \frac{1}{2i\pi} \int_{\Gamma_{\theta,\theta'}} 
    \frac{P}{\lambda (P-\lambda)}d\lambda = \Pi_{\theta,\theta'}(P).
\end{equation}
Hence $\Pi_{\theta,\theta'}(P)$ is a projection.
 
Next, let $\lambda_{0} \in \Sp P$. We may assume that the contour $\Gamma_{(\lambda_{0})}$ does not intersect $\Gamma_{\theta,\theta'}$. Then thanks 
to~(\ref{eq:Appendix-spectral.identity1}) we see that 
$4\pi^{2}\Pi_{\theta,\theta'}(P)\Pi_{\lambda_{0}}(P)$ is equal to
\begin{equation}
    \int_{\Gamma_{\theta,\theta'}}\int_{\Gamma_{(\lambda_{0})}} 
        \frac{P}{\lambda(P-\lambda)(P-\mu)}d\lambda d\mu =  
        \int_{\Gamma_{(\lambda_{0})}} \frac{P}{P-\mu}( 
        \int_{\Gamma_{\theta,\theta'}} \frac{d\lambda}{\lambda(\lambda-\mu)})d\mu.
     \label{eq:Appendix.proof.PilPitt'}
\end{equation}
Therefore, if $\lambda_{0}$ lies outside 
$\Lambda_{\theta,\theta'}$ then $\Pi_{\theta,\theta'}(P)\Pi_{\lambda_{0}}(P)$ is zero, while when $\lambda_{0}$ lies inside $\Lambda_{\theta,\theta'}$ 
using~(\ref{eq:Appendix.proof.trick}) we see that  $\Pi_{\theta,\theta'}(P)\Pi_{\lambda_{0}}(P)$ is equal to
\begin{equation}
     \frac{-1}{2i\pi}\int_{\Gamma_{(\lambda_{0})}} \frac{P}{\mu(P-\mu)}d\mu =  \frac{-1}{2i\pi} \int_{\Gamma_{(\lambda_{0})}} \frac{d\mu}{\mu}¥  +
     \frac{-1}{2i\pi} \int_{\Gamma_{(\lambda_{0})}} 
     \frac{d\mu}{P-\mu}=\Pi_{\lambda_{0}}(P).     
     \label{eq:AppendixC.Pitt'Pil}
\end{equation}
Since $\Pi_{\lambda_{0}}(P)$ has range $E_{\lambda_{0}}(P)$ it then follows that the range of $\Pi_{\theta,\theta'}(P)$ contains
$ E_{\theta,\theta'}(P)$ and its kernel contains $E_{0}(P)\dotplus E_{\theta',\theta+2\pi}(P)$. Hence the result.
\end{proof}

\begin{remark}
    Since Propositions~\ref{prop:Appendix.sectorial-projection-PsiDO} and~\ref{prop:SP.range-kernel} 
    tell us that $\Pi_{\theta,\theta'}(P)$ is a (bounded) \psido\ projection, we see that $\Pi_{\theta,\theta'}(P)$ 
    has either order $0$ or is smoothing.
\end{remark}

\begin{proposition}\label{prop:Appendix-spectral.decomposition-wodzicki-projections}
 Let $L_{\theta_{1}}$ and $L_{\theta_{1}'}$  be spectral cuttings for $P$ and its principal symbol in such way that 
 $\theta' \leq \theta_{1}<\theta'_{1}<\theta+2\pi$. Then:\smallskip

1) The projections $\Pi_{\theta,\theta'}(P)$ and $\Pi_{\theta_{1},\theta_{1}'}(P)$ are disjoint.\smallskip

2) We have $ \Pi_{\theta,\theta'}(P)+\Pi_{\theta',\theta'_{1}}(P)=\Pi_{\theta,\theta'_{1}}(P)$ and $\Pi_{\theta,\theta'}(P)+\Pi_{\theta',\theta+2\pi}(P)=1-\Pi_{0}(P)$. 

\end{proposition}
\begin{proof}
   First, using~(\ref{eq:Appendix-spectral.identity1}) we see that $4\pi^{2}\Pi_{\theta,\theta'}(P)\Pi_{\theta_{1},\theta_{1}'}(P)$ and 
   $4\pi^{2}\Pi_{\theta_{1},\theta_{1}'}(P)\Pi_{\theta_{1},\theta_{1}'}(P)$ 
   are both equal to
\begin{multline}
         \int_{\Gamma_{\theta,\theta'}}\int_{\Gamma_{\theta_{1},\theta_{1}'}} 
       \frac{P^{2}}{\lambda\mu(P-\lambda)(P-\mu)}d\lambda d\mu =\\ 
       \int_{\Gamma_{\theta,\theta'}} \frac{P^{2}}{\lambda(P-\lambda)}
        (\int_{\Gamma_{\theta_{1},\theta_{1}'}}\frac{\mu^{-1}d\mu}{\mu-\lambda})d\lambda        
        + \int_{\Gamma_{\theta_{1},\theta'_{1}}} \frac{P^{2}}{\mu(P-\mu)} 
        (\int_{\Gamma_{\theta,\theta}}\frac{\lambda^{-1}d\mu}{\lambda-\mu})d\mu =0.      
\end{multline}
   Hence $\Pi_{\theta,\theta'}(P)$ and $\Pi_{\theta_{1},\theta_{1}'}(P)$ are disjoint projections.
   
   Next, the operator $\Pi_{\theta,\theta'}(P)+\Pi_{\theta',\theta'_{1}}(P)$ is equal to 
    \begin{equation}
        \frac{1}{2i\pi}\int_{\Gamma_{\theta,\theta'}\cup 
        \Gamma_{\theta',\theta_{1}'}}\frac{P}{\lambda(P-\lambda)}d\lambda 
        = \frac{1}{2i\pi}\int_{\Gamma_{\theta,\theta_{1}'}}\frac{P}{\lambda(P-\lambda)}d\lambda =  \Pi_{\theta,\theta_{1}'}(P),
   \end{equation}
  since integrating $\lambda^{-1}P(P-\lambda)^{-1}$ along $\Gamma_{\theta,\theta'}\cup \Gamma_{\theta',\theta_{1}'}$ is the same as integrating 
   it along $\Gamma_{\theta,\theta_{1}}$. 
In the special case $\theta_{1}'=\theta+2\pi$ the integration along $\Gamma_{\theta,\theta+2\pi}$ reduces to that along the small circle
   $|\lambda|=r$ with clockwise orientation. Therefore, using~(\ref{eq:Appendix.proof.trick}) we see that $ \Pi_{\theta,\theta'}(P)+\Pi_{\theta',\theta+2\pi}(P)$ 
   is equal to
   \begin{equation}
     \frac{1}{2i\pi}\int_{|\lambda|=r}\frac{P}{\lambda(P-\lambda)}d\lambda = 
      \frac{1}{2i\pi}\int_{|\lambda|=r}\frac{d\lambda}{\lambda} +\frac{1}{2i\pi}\int_{|\lambda|=r}\frac{d\lambda}{P-\lambda} =1-\Pi_{0}(P).
   \end{equation}
The proof is thus complete.
\end{proof}

In general, the closures of $ E_{\theta,\theta'}(P)$ and $E_{0}(P)\dotplus E_{\theta',\theta+2\pi}(P)$ don't yield the whole range and the whole kernel of 
$\Pi_{\theta,\theta'}(P)$ but, as we explain in Appendix, there are special cases where they actually do:\smallskip

 
 (i) When the principal symbol of $P$ has no eigenvalues within the angular sector $\theta<\arg \lambda<\theta'$, which is equivalent to 
 $\Pi_{\theta,\theta'}(P)$ being a smoothing operator (see Proposition~\ref{prop:Appendix-CR.smoothingness});\smallskip
 
 (ii) When $P$ is normal, i.e., commutes with its adjoint, and in particular when (see Proposition~\ref{prop:SP.normal});\smallskip
 
 (iii) When $P$ has a complete system of root vectors, that is, the subspace spanned by its root vector is dense (see 
 Proposition~\ref{prop:SP.complete-root-vector-system}).\smallskip  
 
  In the non-normal case it is a difficult issue to determine whether a general closed unbounded operator on a Hilbert space admits a complete system of root 
  vectors. Thanks to a criterion due to Dunford-Schwartz~\cite{DS:LOII} 
  it can be shown that $P$ has a complete system of root 
  vectors when its principal symbol admits spectral cuttings dividing the complex planes into 
  angular sectors of apertures~$<\frac{2 n\pi}{m}$ (see \cite{Ag:EEGEBVP}, \cite{Bu:FPEDO}, \cite{Ag:EOCM}). 
  Therefore, in this case $\Pi_{\theta,\theta'}(P)$ is the projection onto the closure of 
  $ E_{\theta,\theta'}(P)$ and along the closure of $E_{\theta',\theta+2\pi}(P)$.  
  
  In fact, if we content ourselves by determining the range of $\Pi_{\theta,\theta'}(P)$ then it can be shown that the range agrees with the closure
  $E_{\theta,\theta'}(P)$ when we only require the principal symbol of $P$ to admit spectral cuttings dividing the angular sector  $\theta<\arg \lambda<\theta'$ into 
  angular sectors of apertures~$<\frac{2 n\pi}{m}$ (see Proposition~\ref{prop:SP.range-criterion}).

Detailed proofs of the above statements are given in Appendix.

\section{Zeta functions and spectral asymmetry}
\label{sec:SAF}
In this section, we give a detailed review of the spectral asymmetry formulas of Wodzicki~(\cite{Wo:LISA}--\cite{Wo:PhD}) 
for elliptic \psidos. 

Let $P:C^{\infty}(M,\cE) \rightarrow C^{\infty}(M,\cE)$ be an elliptic \psido\ of order $m>0$. Let us first assume that $P$ is selfadjoint. 
Then we have: 
  \begin{equation}
      P_{\scriptscriptstyle{\uparrow}}^{s}=\Pi_{+}(P)|P|^{s}+e^{-i\pi s}\Pi_{-}(P)|P|^{s}, \quad  
      P_{\scriptscriptstyle{\downarrow}}^{s}=\Pi_{+}(P)|P|^{s}+e^{i\pi s}\Pi_{-}(P)|P|^{s}, 
     \label{eq:SAZFNC.PupPdown1}
 \end{equation}
 where $\Pi_{+}(P)$ (resp.~$\Pi_{-}(P)$) denotes the orthogonal projections onto the positive (resp.~negative) eigenspace of 
$P$. Hence we have
 \begin{equation}
       P_{\scriptscriptstyle{\uparrow}}^{s}-P_{\scriptscriptstyle{\downarrow}}^{s} 
      =(e^{-i\pi s}-e^{i\pi s})\Pi_{-}(P)|P|^{s}=(1-e^{2i\pi s})\Pi_{-}(P)P_{\scriptscriptstyle{\uparrow}}^{s}. 
      \label{eq:SAZFNC.PupPdown2} 
 \end{equation}
 Therefore, in the selfadjoint case, the spectral asymmetry of $P$ is encoded by $\Pi_{-}(P)$. 

Suppose now that $P$ is not selfadjoint and let $L_{\theta}=\{\arg \lambda =\theta\}$ and 
$L_{\theta'}=\{\arg \lambda =\theta\}$ be spectral cuttings for both $P$ and its principal 
symbol $p_{m}(x,\xi)$ with $0\leq\theta < \theta'<2\pi$. As observed by Wodzicki~(\cite{Wo:CWSP}, \cite{Wo:PhD}) in this context a substitute to the projection 
$\Pi_{-}(P)$ is provided by the sectorial projection $\Pi_{\theta,\theta'}(P)$ in~(\ref{eq:SP.SAPi}). This stems from:

\begin{proposition}[\cite{Wo:CWSP}, \cite{Wo:PhD}] \label{prop:SAZF.asymmetry-complex-powers}
    For any $s \in \C$ we have 
 \begin{equation}
    P_{\theta}^{s}-P_{\theta'}^{s}=(1-e^{2i\pi s}) \Pi_{\theta,\theta'}(P)P_{\theta}^{s}.
     \label{eq:SAZF.complez-powers.Rtheta-theta'}
\end{equation}
\end{proposition}
\begin{proof}
Since in the integral~(\ref{eq:background.complex-powers-definition}) defining $P_{\theta}^{s}$ 
the value of 
the argument has shifted of $-2\pi$ once $\lambda$ 
has turned around the circle we have
\begin{equation}
    P_{\theta}^{s}=  \frac{e^{2i\pi s}-1}{2i\pi}\int_{\infty}^{r} \frac{x^{s} e^{is(\theta-2\pi)}}{P-xe^{i\theta}} d(xe^{i\theta}) + 
      \int_{\theta}^{\theta-2\pi} \frac{r^{s} e^{ist}}{P-re^{it}}d(re^{it}).
\end{equation}
Similarly, we have 
\begin{equation}
     P_{\theta'}^{s}=  \frac{-e^{2i\pi s}+1}{2i\pi}\int_{r}^{\infty} \frac{x^{s} e^{is(\theta'-2\pi)} }{P-xe^{i\theta'}}d(xe^{i\theta'}) + 
      \int_{\theta'}^{\theta'-2\pi}\frac{ r^{s} e^{ist}}{P-re^{it}}d(re^{it}).     
\end{equation}

Observe that $ \int_{\theta}^{\theta-2\pi} \frac{r^{s} e^{ist}}{P-re^{it}}d(re^{it})- \int_{\theta'}^{\theta'-2\pi} \frac{r^{s} e^{ist}}{P-re^{it}}d(re^{it})$ is equal to
\begin{multline}
    \int_{\theta'-2\pi}^{\theta-2\pi}\frac{ r^{s} e^{ist}}{P-re^{it}}d(re^{it})- \int_{\theta'}^{\theta} \frac{r^{s} e^{ist}}{P-re^{it}}d(re^{it})\\ 
    = (1-e^{2i\pi s}) \int_{\theta'-2\pi}^{\theta-2\pi}\frac{r^{s} e^{ist}}{P-re^{it}}d(re^{it}).  
\end{multline}
Therefore, the operator $P_{\theta}^{s}- P_{\theta'}^{s}$ agrees with
\begin{multline}
 \frac{e^{2i\pi s}-1}{2i\pi}\left(\int_{\infty}^{r} \frac{x^{s} e^{is(\theta-2\pi)}}{P-xe^{i\theta}}d(xe^{i\theta}) 
   +  \int_{\theta-2\pi}^{\theta'-2\pi}\frac{r^{s} e^{ist}}{P-re^{it}}d(re^{it})\right.\\ 
   + \left. \int_{r}^{\infty} \frac{x^{s} e^{is(\theta'-2\pi)}}{P-xe^{i\theta'}}d(xe^{i\theta'})\right).
\end{multline}
In view of the definition~(\ref{eq:SP.SAPi.Gammatt'}) of the contour $\Gamma_{\theta,\theta'}$ this gives
\begin{equation}
     P_{\theta}^{s}-P_{\theta'}^{s}=\frac{e^{2i\pi s}-1}{2i\pi}\int_{\Gamma_{\theta,\theta'}} \lambda^{s}_{\theta}(P-\lambda)^{-1}d\lambda.
    \label{eq:Appendix.proof .asymmetry-formula-intermediate}
\end{equation}

Next, let $\theta_{1}\in (\theta'-2\pi,\theta)$ be such that no eigenvalues of $P$ lie in the sector $\theta_{1}\leq \arg \lambda\leq \theta$. Then 
in the formula~(\ref{eq:background.complex-powers-definition}) for $P_{\theta}^{s}$ we may replace the integration over $\Gamma_{\theta}$ by that over 
a contour $\Gamma_{\theta_{1}}$ defined by~(\ref{eq:background.complex-powers-definition-Gammat})  using $\theta_{1}$ and a radius $r$ smaller than that of 
$\Gamma_{\theta,\theta'}$ in~(\ref{eq:SP.SAPi.Gammatt'}). Thus,
\begin{equation}
    \Pi_{\theta,\theta'}(P)P_{\theta}^{s}=  \frac{1}{4\pi^{2}} \int_{\Gamma_{\theta,\theta'}}\int_{\Gamma_{\theta_{1}}} 
        \lambda^{-1}\mu^{s}_{\theta}P(P-\lambda)^{-1}(P-\mu)^{-1}d\lambda d\mu.
\end{equation}
Using~(\ref{eq:Appendix-spectral.identity1}) we see that $\Pi_{\theta,\theta'}(P)P_{\theta}^{s}$ is equal to
\begin{multline}
\frac{1}{4\pi^{2}}  \int_{\Gamma_{\theta,\theta'}} 
\frac{P}{\lambda (P-\lambda)}(\int_{\Gamma_{\theta_{1}}}\frac{\mu^{s}_{\theta}d\mu}{\mu-\lambda})d\lambda  
        + \frac{1}{4\pi^{2}}\int_{\Gamma_{\theta_{1}}} \frac{\mu^{s}_{\theta} P}{(P-\mu)}(\int_{\Gamma_{\theta,\theta'}}
        \frac{\lambda^{-1}d\lambda}{\lambda-\mu})d\mu\\ 
        = \frac{-1}{2i\pi} \int_{\Gamma_{\theta,\theta'}} \frac{\lambda_{\theta}^{s-1}P}{P-\lambda}d\lambda.  
\end{multline}
Combining this with~(\ref{eq:Appendix.proof.trick}) we obtain:
\begin{equation}
  \label{eq:SAZF.Pi-lambda-Ps}
\Pi_{\theta,\theta'}(P)P_{\theta}^{s}= \frac{-1}{2i\pi}  \int_{\Gamma_{\theta,\theta'}} \lambda_{\theta}^{s-1}d\lambda + 
  \frac{-1}{2i\pi} \int_{\Gamma_{\theta,\theta'}} \!\!
  \frac{ \lambda_{\theta}^{s}}{P-\lambda}d\lambda = \frac{-1}{2i\pi} \int_{\Gamma_{\theta,\theta'}} \!\! \frac{ \lambda_{\theta}^{s}}{P-\lambda}d\lambda.
\end{equation}
Comparing this to~(\ref{eq:Appendix.proof .asymmetry-formula-intermediate}) then gives
\begin{equation}
    P_{\theta}^{s}-P_{\theta'}^{s}=(1-e^{2i\pi s})\Pi_{\theta,\theta'}(P)P_{\theta}^{s}.  
\label{eq:Appendix.asymmetry-complex-powers}
\end{equation}
This proves Proposition~\ref{prop:SAZF.asymmetry-complex-powers} for $\Re s<0$. Since both sides of~(\ref{eq:Appendix.asymmetry-complex-powers}) 
involve holomorphic families of \psidos\ the general case follows by analytic continuation.
\end{proof}

Next,  as the two sides of~(\ref{eq:SAZF.complez-powers.Rtheta-theta'}) are given by holomorphic families of \psidos, from Theorem~\ref{thm:background.TR} 
we immediately get:

\begin{theorem}[\cite{Wo:CWSP}, \cite{Wo:PhD}] \label{thm:SAZF}
We have the equality of meromorphic functions,
    \begin{equation}
     \zeta_{\theta}(P;s)- \zeta_{\theta'}(P;s)=(1-e^{-2i\pi s}) \TR \Pi_{\theta,\theta'}(P)P_{\theta}^{-s}, \quad s \in 
     \C.\label{eq:SAZFNC.SA1} 
   \end{equation}
In particular, at any integer $k \in \Z$ the function $ \zeta_{\theta}(P;s)- \zeta_{\theta'}(P;s)$ is regular 
and there we have 
   \begin{equation}
      \ord P .\lim_{s\rightarrow k} ( \zeta_{\theta}(P;s)- \zeta_{\theta'}(P;s))=2i\pi \Res \Pi_{\theta,\theta'}(P)P^{-k}. 
      \label{eq:SAZFNC.SA2} 
   \end{equation}
\end{theorem}

As a consequence of~(\ref{eq:SAZFNC.SA2}) we see that if at some integer $k$ we have $\Res P^{-k}=0$, so that 
    $\zeta_{\theta}(P;s)$ and $ \zeta_{\theta'}(P;s)$ are regular at 
    $s=k$, then we have: 
    \begin{equation}
       \zeta_{\theta}(P;k)=\zeta_{\theta'}(P;k) \Longleftrightarrow \Res \Pi_{\theta,\theta'}(P)P^{-k}=0.
    \end{equation}

Furthermore, Wodzicki also proved the remarkable result below.

 \begin{theorem}[{\cite[1.24]{Wo:LISA}}]\label{thm:SANCR.independence}
    Let $P:C^{\infty}(M,\cE)\rightarrow C^{\infty}(M,\cE)$ be an elliptic \psido\ of order $m>0$ and let $L_{\theta}=\{\arg \lambda=\theta\}$ be a 
    spectral cutting for $P$ and its principal symbol. Then the regular value $\zeta_{\theta}(P;0)$ is independent of $\theta$. 
  \end{theorem}
\begin{remark}
As with the vanishing of the residue at the origin of the eta function of a selfadjoint elliptic \psido\ Theorem~\ref{thm:SANCR.independence} is not a local result, since 
it is not true that in general the regular value at $s=0$ of the local zeta function $t_{P^{-s}_{\theta}}(x)$ is independent of the spectral cutting 
(see~\cite[pp.~130-131]{Wo:SAZF}). 
\end{remark}

\begin{remark}
The proof of Theorem~\ref{thm:SANCR.independence} in~\cite{Wo:LISA} is quite difficult because it relies on a very involved characterization of local 
invariants of spectral asymmetry. Notice that from~(\ref{eq:SAZFNC.SA2}) we get 
\begin{equation}
     \ord P. ( \zeta_{\theta}(P;0)- \zeta_{\theta'}(P;0))=2i\pi \Res \Pi_{\theta,\theta'}(P),
    \label{eq:SAZFNC.SA3} 
\end{equation}
so that $\zeta_{\theta}(P;0)- \zeta_{\theta'}(P;0)$ is a constant multiple of the noncommutative residue of a \psido\ projection. In fact, 
Wodzcki~\cite[7.12]{Wo:LISA} used Theorem~\ref{thm:SANCR.independence} to prove that the noncommutative residue of a \psido\ projection is always zero. 
However, it follows 
from an observation of Br\"uning-Lesch~\cite[Lem.~2.7]{BL:OEICNLBVP} that the latter result can be deduced in a rather elementary way from the vanishing of the residue 
at the origin of the eta function of a selfadjoint elliptic \psido. Therefore, combining this with~(\ref{eq:SAZFNC.SA3}) 
allows us to prove Theorem~\ref{thm:SANCR.independence} 
without any appeal to Wodzicki's characterization of local invariants of spectral asymmetry. 
\end{remark}

\section{Spectral asymmetry of odd class elliptic \psidos}
\label{sec.odd}
In this section we study the spectral asymmetry of odd class elliptic \psidos. 
Recall that according to~\cite{KV:GDEO} a \psido\ $Q$ of integer order $m$ is an odd class \psido\ when, in local 
trivializing coordinates, its symbol $q(x,\xi)\sim \sum_{j\geq 0} q_{m-j}(x,\xi)$ is polyhomogeneous with respect to the dilation by $-1$, i.e., for $j=0,1,\ldots$ we have 
\begin{equation}
    q_{m-j}(x,-\xi)=(-1)^{m-j}q_{m-j}(x,\xi).
     \label{eq:SAZF.transmission-property}
\end{equation}
This gives rise to a subalgebra of $\Psi^{\Z}(M,\cE)$ which contains 
all the differential operators and the parametrices of elliptic odd class \psidos. 

Moreover, the condition 
$q_{-n}(x,-\xi)=(-1)^{n}q_{-n}(x,\xi)$  implies that, 
when the dimension of $M$ is odd, the noncommutative residue of an odd class \psido\ vanishes locally, i.e., the 
density $c_{Q}(x)$ given by~(\ref{eq:NCR.density}) vanishes.

 \begin{theorem}\label{SAZF.odd-PsiDOs1}
    Suppose that $\dim M$ is odd and that $P$ is an odd class \psido\ of even integer order $m\geq 2$.  Then 
    $\zeta_{\theta}(P;s)$ is regular at every integer point and its values there are independent of the cutting.
  \end{theorem}
\begin{proof}
 In some local trivializing  coordinates let $p(x,\xi)\sim \sum_{j\geq 0}p_{m-j}(x,\xi)$ denote the symbol of $P$ and let 
  $q(x,\xi,\lambda)\sim \sum_{j\geq 0} q_{-m-j}(x,\xi,\lambda)$ be  the symbol with parameter of $(P-\lambda)^{-1}$ as in~\cite{Se:CPEO}, so that 
  $q_{-m-j}(x,t\xi,t^{m}\lambda)=t^{-m-j}q_{-m-j}(x,\xi,\lambda)$ for $t\neq 0$ and $\sim$ is taken in the sense of symbols with parameter 
  of~\cite[p.~295]{Se:CPEO}.  
  Then by~(\ref{eq:Appendix.symbol-Wodzicki-projection}) 
  the symbol $\pi(x,\xi) \sim \sum_{j\geq 0} \pi_{-j}(x,\xi)$ is given by 
 \begin{equation}
     \pi_{-j}(x,\xi)= \frac{-1}{2i\pi} \int_{\Gamma_{(x,\xi)}} q_{-m-j}(x,\xi;\lambda) d\lambda,
     \label{eq:SAZF.symbol-Wodzicki-projection}
\end{equation}
 where $\Gamma_{(x,\xi)}$ is a direct-oriented bounded contour contained in the angular sector $\Lambda_{\theta,\theta'}=\{\theta < \arg \lambda<\theta'\}$ 
 which isolates from $\C\setminus \Lambda_{\theta,\theta'}$ the eigenvalues of $p_{m}(x,\xi)$ that lie in 
 $\Lambda_{\theta,\theta'}$. 
 
 At the level of symbols the equality $(P-\lambda)(P-\lambda)^{-1}=1$ gives
\begin{equation}
   1=  p\# (q-\lambda) \sim (p(x,\xi)-\lambda)q(x,\xi,\lambda)+ \sum_{\alpha\neq 0}\frac{1}{\alpha!}
     \partial_{\xi}^{\alpha}p(x,\xi)D_{x}^{\alpha}q(x,\xi,\lambda).
\end{equation}
From this we get
  \begin{equation}
      q_{-m}(x,\xi, \lambda)=(p_{m}(x,\xi)-\lambda)^{-1},\label{ref:SAZF.symbol-resolvent1}
  \end{equation}
and for $j=1,2,\ldots$ we see that $q_{-m-j}(x,\xi, \lambda)$ is equal to
  \begin{equation}
      -(p_{m}(x,\xi)-\lambda)^{-1}\!\! \sum_{\substack{|\alpha|+k+l=j,\\ l\neq 
    j}}\!\! \frac{1}{\alpha!} \partial_{\xi}^{\alpha}p_{m-k}(x,\xi)D_{x}^{\alpha}q_{-m-l}(x,\xi,\lambda).
    \label{ref:SAZF.symbol-resolvent2}
 \end{equation}
Since the symbol $p(x,\xi)$ satisfies~(\ref{eq:SAZF.transmission-property}), it follows 
from~(\ref{ref:SAZF.symbol-resolvent1}) and  (\ref{ref:SAZF.symbol-resolvent2}) that 
for $j=0,1,\ldots$ we have
 \begin{equation}
     q_{-m-j}(x,-\xi, (-1)^{m}\lambda)=(-1)^{-m-j}q_{-m-j}(x,\xi,\lambda).
        \label{eq:SAZF.transmission-property-resolvent}
 \end{equation}

Now, assume $n$ is odd and $m$ is even. As alluded to above the noncommutative residue of an odd class \psido\ 
is zero in odd dimension. Since the odd class \psidos\ form an algebra containing all the parametrices of odd class elliptic \psidos\  
it follows that for any integer $k$ the operator 
$P^{-k}$ is an odd class \psido\ and its noncommutative residue is zero. Therefore, the zeta functions $\zeta_{\theta}(P;s)$ 
and $\zeta_{\theta'}(P;s)$ are regular at all integer points. 

On the other hand, since $m$ is even thanks to~(\ref{eq:SAZF.transmission-property-resolvent}) we see that
\begin{multline}
   \pi_{-j}(x,-\xi)= \frac{-1}{2i\pi}\int_{\Gamma_{(x,-\xi)}} \!\! q_{-m-j}(x,-\xi, \lambda)d\lambda=\\
        \frac{(-1)^{m-j+1}}{2i\pi}\int_{\Gamma_{(x,\xi)}}\!\!  q_{-m-j}(x,\xi, \lambda)d\lambda =(-1)^{-j} \pi_{-j}(x,\xi).
\end{multline}
Hence $\Pi_{\theta,\theta'}(P)$ is an odd class \psido.  Therefore, for any $k\in \Z$ the operator 
$\Pi_{\theta,\theta'}(P)P^{-k}$ is an odd class \psido\ as well, and so $\Res \Pi_{\theta,\theta'}(P)P^{-k}=0$. It then   
follows from Theorem~\ref{thm:SAZF} that $\zeta_{\theta}(P;k)=\zeta_{\theta'}(P;k)$.  
\end{proof}

 \begin{theorem}\label{SAZF.odd-PsiDOs2}
    Assume $\dim M$ is even, $P$ is an odd class \psido\ of odd integer order $m\geq 1$ such  all the eigenvalues of its principal symbol 
      lie in the open cone $\{\theta <\arg \lambda<\theta'\}\cup\{\theta +\pi<\arg \lambda<\theta'+\pi\}$.  Then:\smallskip  
      
     1) For any integer $k \in \Z$ we have  
      \begin{equation}
         \ord P. \lim_{s\rightarrow k}(\zeta_{\theta}(P;s)- \zeta_{\theta'}(P;s))=i\pi  \Res P^{-k}.
          \label{eq:SAZF.odd-PsiDOs.m-odd}
      \end{equation}
      
     2) At every integer at which they are not singular the  functions 
      $\zeta_{\theta}(P;s)$ and $\zeta_{\theta'}(P;s)$ take on 
       the same regular value.
 \end{theorem}
\begin{proof}
Since all the eigenvalues of $p_{m}(x,\xi)$ are contained in 
the cone $\mathcal{C}_{\theta,\theta'}:=\{\theta <\arg \lambda<\theta'\}\cup\{\theta +\pi<\arg \lambda<\theta'+\pi\}$. Then $P$ has at most finitely 
many eigenvalues in $\mathcal{C}_{\theta,\theta'}$ and by Proposition~\ref{prop:Appendix-CR.smoothingness} 
the sectorial projections $\Pi_{\theta',\theta+\pi}(P)$ and $\Pi_{\theta'+\pi,\theta+2\pi}(P)$ are smoothing operators.  

On the other hand, by Proposition~\ref{prop:Appendix-spectral.decomposition-wodzicki-projections} we have
\begin{equation}
     \Pi_{\theta,\theta'}(P)+\Pi_{\theta',\theta+\pi}(P)+\Pi_{\theta+\pi,\theta'+\pi}(P)+\Pi_{\theta'+\pi,\theta+2\pi}(P)=1-\Pi_{0}(P).
     \label{eq:SAZF.decomposing-Wodzicki-projections}
\end{equation}
Since $\Pi_{\theta',\theta+\pi}(P)$ and $\Pi_{\theta'+\pi,\theta+2\pi}(P)$, as well as $\Pi_{0}(P)$, are smoothing operators it follows that 
\begin{equation}
    \Pi_{\theta,\theta'}(P)+\Pi_{\theta+\pi,\theta'+\pi}(P)=1 \quad \bmod \Psi^{-\infty}(M,\cE).
\end{equation}
Combining this with~(\ref{eq:SAZF.symbol-Wodzicki-projection}) we see that at the level of symbols we get
     \begin{gather}
         \frac{-1}{2i\pi} \int_{\Gamma_{(x,\xi)}} q_{-m}(x,\xi,\lambda)d\lambda + \frac{-1}{2i\pi} 
         \int_{-\Gamma_{(x,\xi)}}q_{-m}(x,\xi,\lambda)d\lambda =1, \label{eq:SAZF.m-odd.n-even1}\\
         \frac{-1}{2i\pi} \int_{\Gamma_{(x,\xi)}} q_{-m-j}(x,\xi,\lambda)d\lambda + \frac{-1}{2i\pi} 
         \int_{-\Gamma_{(x,\xi)}}q_{-m-j}(x,\xi,\lambda)d\lambda =0, \quad j\geq 
         1.\label{eq:SAZF.m-odd.n-even2}
     \end{gather}

Next, observe that the formula~(\ref{eq:SAZF.transmission-property-resolvent}) in the proof of Theorem~\ref{SAZF.odd-PsiDOs1} 
is actually true independently of the parities of $m$ and $n$.  Therefore, we may combine it with~(\ref{eq:SAZF.m-odd.n-even1}) to get
\begin{multline}
  \pi_{0}(x,-\xi) - 1= \frac{1}{2i\pi}\int_{-\Gamma_{(x,\xi)}}\!\!q_{-m}(x,-\xi,\lambda)d\lambda  \\
         = \frac{-1}{2i\pi}\int_{\Gamma_{(x,\xi)}}q_{-m}(x,-\xi,-\lambda)d\lambda = (-1)^{m} \pi_{0}(x,\xi)= 
-\pi_{0}(x,\xi).
\end{multline}

Similarly, using~(\ref{eq:SAZF.transmission-property-resolvent})  and (\ref{eq:SAZF.m-odd.n-even2}) for $j=1,2,\ldots$ we get  
\begin{multline}
    \pi_{-j}(x,-\xi) = \frac{-1}{2i\pi} \int_{-\Gamma_{(x,\xi)}}q_{-m-j}(x,-\xi,-\lambda)d\lambda \\
   = (-1)^{m-j}\pi_{-j}(x,\xi)= (-1)^{j+1}\pi_{-j}(x,\xi).   
\end{multline}

Now, let $k \in \Z$ and let $p^{(k)}\sim \sum_{j\geq 0} p^{(k)}_{-km-j}$ denote the symbol of $P^{-k}$. Then 
the symbol $r^{(k)}_{-n}$ of degree $-n$ of $R^{(k)}=\Pi_{\theta,\theta'}(P)P^{-k}$ is given by 
\begin{equation}
    r^{(k)}_{-n}(x,-\xi)=\!\!\sum_{|\alpha|+j+l=n-km}\!\! \frac{1}{\alpha !} 
\partial_{\xi}^{\alpha}\pi_{-j}(x,\xi)D_{x}^{\alpha}p^{(k)}_{-km-l}(x,\xi).
     \label{eq:SAZF.symbol-r(k)(-n)}
\end{equation}
Since $P^{-k}$ is an odd class \psido, using~(\ref{eq:SAZF.m-odd.n-even1}) and  (\ref{eq:SAZF.m-odd.n-even2})   we obtain:
\begin{multline}
      r^{(k)}_{-n}(x,-\xi) = \sum_{j+l+|\alpha|=n-km}\!\! \frac{1}{\alpha !} 
(\partial_{\xi}^{\alpha}\pi_{-j})(x,-\xi)
    (D_{x}^{\alpha}p^{(k)}_{-km-l})(x,-\xi)\\
= \sum_{l+|\alpha|=n-km}\!\! \frac{(-1)^{|\alpha|-km-l}}{\alpha !} 
\partial_{\xi}^{\alpha}[1-\pi_{0}(x,\xi)]    D_{x}^{\alpha}p^{(k)}_{-km-l})(x,\xi)\\
  -\sum_{j+l+|\alpha|=n-km}\!\! \frac{(-1)^{j+|\alpha|-km-l}}{\alpha !} 
\partial_{\xi}^{\alpha}\pi_{-j}(x,\xi) D_{x}^{\alpha}p^{(k)}_{-km-l}(x,\xi)\\ 
   =  (-1)^{n}p_{-n}^{(k)}(x,\xi) -(-1)^{n} \!\!\sum_{|\alpha|+j+l=n-km} \!\! 
    \frac{1}{\alpha !} (\partial_{\xi}^{\alpha}\pi_{-j})(x,\xi) 
    (D_{x}^{\alpha}p^{(k)}_{-km-l})(x,\xi).
\end{multline}
Combining this with~(\ref{eq:SAZF.symbol-r(k)(-n)}) and the fact that 
$n$ is even we get 
\begin{equation}
    r^{(k)}_{-n}(x,\xi)+r^{(k)}_{-n}(x,-\xi)=p_{-n}^{(k)}(x,\xi).
     \label{eq:SAZF.rkxi-rk-xi}
\end{equation}  

Moreover, we have 
\begin{equation}
    \int_{|\xi|=1}r^{(k)}_{-n}(x,-\xi)d^{n-1}\xi=(-1)^{n}\int_{|\xi|=1}r^{(k)}_{-n}(x,\xi)d^{n-1}\xi 
        =(2\pi)^{-n}c_{R^{(k)}}(x), 
\end{equation}
where $c_{R^{(k)}}(x)$ is the residual density~(\ref{eq:NCR.density}). Thus~(\ref{eq:SAZF.rkxi-rk-xi}) yields $2c_{R^{(k)}}(x)=c_{P^{-k}}(x)$, 
from which we get $\Res \Pi_{\theta,\theta'}(P)P^{-k}=\Res R^{(k)}=\frac{1}{2} \Res P^{-k}$. Combining this with 
Theorem~\ref{thm:SAZF} then gives
\begin{equation}
    \ord P.\lim_{s\rightarrow k}(\zeta_{\theta}(P;s)- \zeta_{\theta'}(P;s))=2i\pi \Res P^{-k}. 
\end{equation}

Finally, by Proposition~\ref{prop:background.zeta-function} the functions $\zeta_{\theta}(P;s)$ and 
$\zeta_{\theta'}(P;s)$ are 
regular at $k\in\Z$ iff $\Res P^{-k}=0$. As $\ord P.\lim_{s\rightarrow 
k}(\zeta_{\theta}(P;s)- \zeta_{\theta'}(P;s))=2i\pi \Res P^{-k}$  it follows that whenever 
$\zeta_{\theta}(P;s)$ and $\zeta_{\theta'}(P;s)$ are regular at an integer their regular values there coincide. 
In particular, as they are always regular at the origin we have $\zeta_{\theta}(P;0)= \zeta_{\theta'}(P;0)$. 
\end{proof}

\begin{remark}\label{rem:SAZF.diff-op}
    As the noncommutative residue of a differential operator is always zero, we see that  if in  
   Theorems~\ref{SAZF.odd-PsiDOs1} and~\ref{SAZF.odd-PsiDOs2} 
   we further assume that $P$ is a differential operator, then at every integer not between $1$ and $\frac{n}{m}$ the  functions 
  $\zeta_{\theta}(P;s)$ and $\zeta_{\theta'}(P;s)$ are non-singular and share the same regular value. 
\end{remark}

Finally, the proofs of Theorems~\ref{SAZF.odd-PsiDOs1} and~\ref{SAZF.odd-PsiDOs2} are based on the analysis of the symbol of 
$\Pi_{\theta,\theta'}(P)$, so the theorems 
ultimately hold at level the local zeta functions 
$\zeta_{\theta}(P;s)(x):=\tr_{\cE}t_{P_{\theta}^{-s}}(x)$. In particular, for the regular value at $s=0$ we get: 

\begin{theorem}\label{SAZF.odd-PsiDOs-local}
   If $P$ satisfies either the assumptions of Theorem~\ref{SAZF.odd-PsiDOs1} or that of Theorem~\ref{SAZF.odd-PsiDOs2}, 
 then $\zeta_{\theta}(P;0)(x)$ is independent of the cutting. 
\end{theorem}

 This shows that the independence of  $\zeta_{\theta}(P;0)(x)$ with respect to the cutting, while not true in general (see~\cite[pp.~130-131]{Wo:SAZF}), 
 nevertheless occurs for a large  class of elliptic \psidos. 
%

\section{Spectral asymmetry of selfadjoint odd class elliptic \psidos}
\label{sec.selfadjoint-odd}
In this section we specialize the results from the previous sections to selfadjoint odd class elliptic \psidos\ and use them to study the eta function of such 
operators.

Let $P:C^{\infty}(M,\cE)\rightarrow 
C^{\infty}(M,\cE)$ be a selfadjoint odd class elliptic \psido\ of integer order $m\geq 1$. 
Since the principal symbol $p_{m}(x,\xi)$ of $P$ is selfadjoint, the assumption in Theorem~\ref{SAZF.odd-PsiDOs2} 
on the location of the eigenvalues of $p_{m}$ is 
always satisfied if we take $0<\theta<\pi<\theta'<2\pi$. 

Now, Theorems~\ref{SAZF.odd-PsiDOs1} and~\ref{SAZF.odd-PsiDOs2} tell us
that if $\dim M$ and $\ord P$ have opposite parities then there are many integer points at which the zeta functions $\zetaup(P;s)$ and 
$\zetadown(P;s)$ are not asymmetric. However, they also  allow us to single out points at which the 
asymmetry of zeta functions always occurs. For instance, we have: 

\begin{proposition}\label{prop:SAZF.selfadjoint-asymmetry}
   If $\dim M$ is even and $P$ is an odd class selfadjoint elliptic \psido\ of order 1, then 
   we always have $ \lim_{s\rightarrow n}\frac{1}{i}(\zetaup(P;s)- \zetadown(P;s))> 0$.
\end{proposition}
\begin{proof}
By Theorem~\ref{SAZF.odd-PsiDOs2}  we have $\lim_{s\rightarrow n}\frac{1}{i}(\zetaup(P;s)- \zetadown(P;s))=\pi \Res 
P^{-n}$. 
Moreover, since  $P^{-n}$ has order $-n$ its symbol of degree $-n$ is its principal symbol $p_{m}(x,\xi)^{-n}$, so we have 
$\Res P^{-n}=(2\pi)^{-n}\int_{S^{*}M}\tr p_{m}(x,\xi)^{-n}dxd\xi$,
where $S^{*}M$ denotes the cosphere bundle of $M$ with its induced metric. 

On the other hand, as $p_{m}(x,\xi)$ is selfadjoint 
and $n$ is even we have 
$\tr p_{m}(x,\xi)^{-n}=\tr[p_{m}(x,\xi)^{-\frac{n}{2}*}p_{m}(x,\xi)^{-\frac{n}{2}}]>0$. Hence
$\Res P^{-n}$ and $\lim_{s\rightarrow n}\frac{1}{i}(\zetaup(P;s)- \zetadown(P;s))$ are positive numbers.
\end{proof}

Next, as observed by Shubin~\cite[p.~114]{Sh:POST} (see also \cite[p.~116]{Wo:SAZF}), we can relate $\zetaup(P;s)- \zetadown(P;s)$ to the 
eta function $\eta(P;s)$ as follows. Let $F=\Pi_{+}(P)-\Pi_{-}(P)$ be the sign operator of $P$. Then using~(\ref{eq:SAZFNC.PupPdown1}) we get:
\begin{equation}
   P_{\scriptscriptstyle{\uparrow}}^{s}-F|P|^{s}=(1+e^{-i\pi s})\Pi_{-}(P)|P|^{s}.
\end{equation}
Combining this with~(\ref{eq:SAZFNC.PupPdown2}) and the fact that 
$(1-e^{i\pi s})(e^{-i\pi s}+1)=e^{-i\pi s}-e^{i\pi s}$ we obtain
\begin{equation}
    P_{\scriptscriptstyle{\uparrow}}^{s}-P_{\scriptscriptstyle{\downarrow}}^{s} = (e^{-i\pi s}-e^{i\pi s}) (1+e^{-i\pi 
    s})^{-1}(P_{\scriptscriptstyle{\uparrow}}^{s}-F|P|^{s})=
    (1-e^{i\pi s})( 
    P_{\scriptscriptstyle{\uparrow}}^{s}-F|P|^{s}).  
\end{equation}
Since $\eta(P;s)=\TR F|P|^{-s}$ we get:

\begin{proposition}\label{prop:SAZF.selfadjoint-case-eta}
1) We have the equality of meromorphic functions,
  \begin{equation}
     \zeta_{\scriptscriptstyle{\uparrow}}(P;s)-\zeta_{\scriptscriptstyle{\downarrow}}(P;s)= (1-e^{-i\pi s}) 
     \zeta_{\scriptscriptstyle{\uparrow}}(P;s) - (1-e^{-i\pi s})\eta(P;s), \quad s\in \C.
     \label{eq:Odd.zeta-eta}
  \end{equation} 
In particular, for any $k \in \Z$ we have
  \begin{equation}
     \ord P. \lim_{s\rightarrow k}(\zeta_{\scriptscriptstyle{\uparrow}}(P;s)-\zeta_{\scriptscriptstyle{\downarrow}}(P;s)) = 
i\pi \Res P^{-k}- i\pi \ord P.\res_{s=k}\eta(P;s). 
     \label{eq:Odd.eta-zeta-up-down}
 \end{equation} 
 
2) Let $k\in \Z$ and suppose that $\Res P^{-k}=0$, so that 
$\zeta_{\scriptscriptstyle{\uparrow}}(s)$ and $\zeta_{\scriptscriptstyle{\downarrow}}(s)$ 
are both regular at $s=k$. Then we have:
    \begin{equation}
        \zeta_{\scriptscriptstyle{\uparrow}}(P;k)=\zeta_{\scriptscriptstyle{\downarrow}}(P;k) \Longleftrightarrow 
        \text{$\eta(P;s)$ is regular at $s=k$}.
        \label{eq:SA.zeta-eta3}
     \end{equation}
\end{proposition}

Now, by a well known result of Branson-Gilkey~\cite{BG:REFODT} in even dimension the eta function of a geometric Dirac 
operator is an entire function. In fact, the latter is a special case of the more general  result below. 

\begin{theorem}\label{thm:SAZF.regularity-eta}
    1) If  $\dim M$ and $\ord P$ have opposite parities then
    $\eta(P;s)$ is regular at every integer point. \smallskip 
    
    2) If $P$ has order 1 and  $\dim M$ is even then $\eta(P;s)$ is an entire function.
\end{theorem}
\begin{proof}
    Let $ k\in \Z$. Since $\dim M$ and $\ord P$ have opposite parities Theorem~\ref{SAZF.odd-PsiDOs1} and Theorem~\ref{SAZF.odd-PsiDOs2} tell us that
    $i\pi \Res P^{-k}$ and $\ord P.\lim_{s\rightarrow  k} 
        (\zeta_{\scriptscriptstyle{\uparrow}}(P;s)-\zeta_{\scriptscriptstyle{\downarrow}}(P;s))$ in~(\ref{eq:Odd.eta-zeta-up-down}) 
        either are both equal to zero (when $\dim M$ is odd and $\ord P$ is even) or are equal to each other (when $\dim M$ is even and 
        $\ord P$ is odd). In any case~(\ref{eq:Odd.eta-zeta-up-down}) shows that $\eta(P;s)$ is regular at $s=k$. 
        
  On the other hand, when $P$ has order 1 Proposition~\ref{prop:background.eta-function} implies that $\eta(P;s)$ 
is holomorphic on $\C\setminus \Z$. Thus, when $\dim M$ is even and $P$ has order $1$ the function $\eta(P;s)$ is entire. 
\end{proof}

\begin{remark}
Theorem~\ref{thm:SAZF.regularity-eta} has  been obtained independently by Grubb~\cite{Gr:RATZLE} using a different approach. 
\end{remark}

\begin{remark}\label{rem:SAZF.Bruening-Seeley}
Theorem~\ref{thm:SAZF.regularity-eta}  allows us to simplify in the odd dimensional case the index formula of Br\"uning-Seeley~\cite[Thm.~4.1]{BS:ITFORSO} 
for a first order elliptic differential operator on a manifold $M$ with cone-like singularities. The contribution of the singularities to this formula involves 
the residues at integer  points of some first order selfadjoint elliptic differential operators on manifolds of dimension $\dim M-1$. Thus when $\dim M$ 
is odd Theorem~\ref{thm:SAZF.regularity-eta} insures us that all these residues are zero, hence disappear from the formula. 
\end{remark}

\section{Spectral asymmetry of Dirac operators}
\label{sec.Dirac}
In this section we make use of the results of the previous sections to express in geometric terms the spectral asymmetry of Dirac operators. 

Throughout all the section we assume that $\dim M$ is even and that $\cE$ is endowed with a Clifford module structure, that is, 
a $\Z_{2}$-grading $\cE=\cE_{+}\oplus \cE_{-}$ and an action of the Clifford bundle $\Cl(M)$ on $\cE$ anticommuting with each other. 
Given a unitary Clifford connection $\nabla^{\cE}$ we get a Dirac operator $\sD_{\cE}$ as the composition, 
\begin{equation}
    \sD_{\cE}: C^{\infty}(M,\cE) \stackrel{\nabla^{\cE}}{\longrightarrow} C^{\infty}(M,T^{*}M\otimes \cE) 
    \stackrel{c}{\longrightarrow} C^{\infty}(M,\cE),
\end{equation}
where $c$ denotes the Clifford action of $T^{*}M$ on $\cE$ (see~\cite[Sect.~3.3]{BGV:HKDO}). 

This setting covers many geometric examples, e.g., 
the Dirac operator on a spin Riemannian manifold withcoefficients in a Hermitian vector 
bundle,  the Gauss-Bonnet and signature operators on an oriented Riemannian manifold, or even the $\overline{\partial}+\overline{\partial}^{*}$-operator on a 
Kaehler manifold. 

The main result of this section is the following.

 \begin{theorem}\label{thm:SAZF.Dirac-operator}
     1) The function $\zetaup(\sD_{\cE};s)- \zetadown(\sD_{\cE};s)$ is entire.\smallskip
    
    2)  At every odd integer and at every even integer not 
     between $2$ and $n$ the functions $\zetaup(\sD_{\cE};s)$ and $\zetadown(\sD_{\cE};s)$ are regular and have the same regular value.\smallskip 
     
     3) For $k=2,4,\ldots,n$ we have 
     \begin{equation}
         \lim_{s\rightarrow k}(\zetaup(\sD_{\cE};s)- \zetadown(\sD_{\cE};s))=i\pi 
                  \int_{M}A_{k}(R^{M},F^{\cE/\sS})(x) \sqrt{g(x)}d^{n}x,
          \label{eq:SAZF.Dirac-operator.universal-formula}
     \end{equation}
     where $A_{k}(R^{M},F^{\cE/\sS})(x)$ is a universal polynomial in complete tensorial contractions of the covariant derivatives of the 
    Riemannian curvature $R^{M}$ of $M$ and of the twisted curvature $F^{\cE/\sS}$ of $\cE$ as defined in~\cite[Prop.~3.43]{BGV:HKDO}. 
     In particular, 
     \begin{gather}
         \lim_{s\rightarrow n}(\zetaup(\sD_{\cE};s)-\zetadown(\sD_{\cE};s))= 
         2i\pi (4\pi)^{-n/2}\Gamma(\frac{n}{2})^{-1}\op{rk} \cE.\op{vol}M,  \label{eq:SAZF.Dirac-operator.n}\\  
          \lim_{s\rightarrow n-2}(\zetaup(\sD_{\cE};s)- \zetadown(\sD_{\cE};s))= -i\pi c_{n} \op{rk} \cE
          \int_{M}r_{M}(x) \sqrt{g(x)}d^{n}x,
          \label{eq:SAZF.Dirac-operator.n-2}
     \end{gather}
      where $c_{n}=\frac{1}{12} (n-2)(4\pi)^{-n/2}\Gamma(\frac{n}{2})^{-1}$ and $r_{M}$ denotes the scalar curvature of $M$.
 \end{theorem}
 \begin{proof}
     First, as $\sD_{\cE}$ is a first order differential operator 
     Proposition~\ref{prop:background.zeta-function} tells us that the function $\zetaup(\sD_{\cE};s)- \zetadown(\sD_{\cE};s)$ can have poles only at integer 
     points and by Theorem~\ref{thm:SAZF} the function is regular at these points. Thus $\zetaup(\sD_{\cE};s)- \zetadown(\sD_{\cE};s)$ is an entire 
     function. 
     
     Second, since $n$ is even it follows from  
     Theorem~\ref{SAZF.odd-PsiDOs2} and Remark~\ref{rem:SAZF.diff-op}  that at every  integer $k$ not between $1$ and 
     $n$ the  functions   $\zetaup(\sD_{\cE};s)$ and  $\zetadown(\sD_{\cE};s)$ are regular and have the same regular value.
     
     Next, by construction $\sD_{\cE}$ anticommutes with the $\Z_{2}$-grading of $\cE$, so when $k$ is odd $\sD_{\cE}^{-k}$ also anticommutes with the 
     $\Z_{2}$-grading. At the level of the residual density $c_{\sD^{-k}_{\cE}}(x)$ this  implies that it 
     take values in  endomorphisms of $\cE$ intertwining $\cE^{+}$ and $\cE^{-}$, so that we have 
     $\tr_{\cE}c_{\sD^{-k}_{\cE}}(x)=0$ and $\Res P^{-k}$ vanishes. Thus, at $s=k$ the functions $\zetaup(\sD_{\cE};s)$ and  
     $\zetadown(\sD_{\cE};s)$ are regular and so have same regular value by Theorem~\ref{SAZF.odd-PsiDOs2}.
     
     Now, let us assume that $k=2l$ for some integer $l$ between $0$ and $\frac{n}{2}$. Thanks to~(\ref{eq:SAZF.odd-PsiDOs.m-odd}) we have 
     \begin{equation}
          \lim_{s\rightarrow k}(\zetaup(\sD_{\cE};s)- \zetadown(\sD_{\cE};s))=i\pi \int_{M}\tr_{\cE}c_{(\sD_{\cE}^{2})^{-l}}(x). 
          \label{eq:Dirac.zetaup-zetadown.2l}
     \end{equation}

As it is well-known (see~\cite[3.23]{Wo:NCRF}) the densities $c_{(\sD_{\cE}^{2})^{-l}}(x)$, 
     $l=1,\ldots,\frac{n}{2}$, are related to the coefficient of the small time heat-kernel asymptotics, 
     \begin{equation}
         k_{t}(x,x)\sim t^{-\frac{n}{2}} \sum_{j\geq 0} t^{j}a_{j}(\sD_{\cE}^{2})(x) \qquad \text{as $t\rightarrow 0^{+}$},
         \label{eq:SAZF.heat-kernel-asymptotics}
     \end{equation}
     where $k_{t}(x,y)$, $t>0$, denotes the heat kernel of $\sD_{\cE}^{2}$. More precisely, we  have 
      \begin{equation}
          c_{(\sD_{\cE}^{2})^{-l}}(x) = \frac{2}{(l-1)!} a_{\frac{n}{2}-l}(\sD_{\cE}^{2})(x). 
           \label{eq:SAZF.NCR-heat-kernel-asymptotics}
     \end{equation}
    
     On the other hand, the operator $\sD_{\cE}^{2}$ is a Laplace type operator, since by the Lichnerowicz's formula
     we have $ \sD_{\cE}^{2}=(\nabla^{\cE})^{*}\nabla^{\cE}+c(F^{\cE/\sS})+\frac{1}{4}r_{M}$ (see~\cite[Thm.~3.52]{BGV:HKDO}). Therefore, 
     by~\cite[pp.~334-336]{Gi:ITHEASIT} 
 each density $a_{j}(\sD_{\cE}^{2})(x)$'s is of the form 
     $\tilde{A}_{j}(R^{M},F^{\cE/\sS})\sqrt{g(x)}dx$, for some universal polynomial $\tilde{A}_{j}(R^{M},F^{\cE/\sS})$ 
     in complete tensorial contractions of the curvatures $R^{M}$ and $F^{\cE/\sS}$. In particular, we have 
     \begin{gather}
        \tilde{A}_{0}(R^{M},F^{\cE/\sS})= (4\pi)^{-n/2}\op{id}_{\cE},\\
         \tilde{A}_{1}(R^{M},F^{\cE/\sS})=\frac{-(4\pi)^{-n/2}}{12}(r_{M}\op{id}_{\cE}+2c(F^{\cE/\sS})).
     \end{gather}
    Combining this with~(\ref{eq:Dirac.zetaup-zetadown.2l}) and~(\ref{eq:SAZF.NCR-heat-kernel-asymptotics}) and the fact that $ 
    \Tr_{\cE}c(F^{\cE/\sS})=0$ then gives the 
    formulas~(\ref{eq:SAZF.Dirac-operator.universal-formula})--(\ref{eq:SAZF.Dirac-operator.n-2}). 
\end{proof}

As an immediate consequence of~(\ref{eq:SAZF.Dirac-operator.n}) we get: 

\begin{proposition}\label{prop:SAZF.Dirac-asymmetry} 
1) The value of $\lim_{s\rightarrow n-2}(\zetaup(\sD_{\cE};s)- \zetadown(\sD_{\cE};s))$ is independent of the 
Clifford data $(\cE,\nabla^{\cE})$.\smallskip 

2) If $ \int_{M}r_{M} \sqrt{g(x)}dx \neq 0$ then we have
 $\lim_{s\rightarrow n-2}(\zetaup(\sD_{\cE};s)- \zetadown(\sD_{\cE};s)) \neq 0$ for any Clifford data $(\cE,\nabla^{\cE})$.
\end{proposition}

Finally, the integral $\int_{M}r_{M}(x)\sqrt{g(x)}dx$ is the 
    Einstein-Hilbert action of the metric $g$,  which gives the contribution of gravity forces to the action functional in
    general relativity. Therefore, it is an important issue in noncommutative geometry and mathematical physics to give an operator theoretic 
    formulation of this action. The first one by given by Connes~\cite{Co:GCMFNCG} in terms 
    of $\Res \sD_{M}^{-n+2}$ (see also~\cite{KW:GNGWR}, \cite{Ka:DOG}), 
    but we see here that thanks to~(\ref{eq:SAZF.Dirac-operator.n}) we get another spectral interpretation of the Einstein-Hilbert action.   

    \appendix 

    \section*{Appendix}
 \setcounter{section}{1}
 In this appendix we gather the main results regarding the spectral interpretation of the sectorial projection of an elliptic \psido.
 
Let $P:C^{\infty}(M,\cE)\rightarrow 
C^{\infty}(M,\cE)$ be an elliptic \psido\ of order $m>0$ 
and assume that $L_{\theta}=\{\arg \lambda =\theta\}$  and 
$L_{\theta'}=\{\arg \lambda =\theta\}$ are spectral cuttings for both $P$ and its principal 
symbol $p_{m}(x,\xi)$ with $\theta< \theta'\leq \theta+2\pi$. We let $\Pi_{\theta,\theta'}(P)$ be the corresponding sectorial projection as defined 
in~(\ref{eq:SP.SAPi}) and we shall use in the sequel the notation introduced in Section~\ref{sec.SP}. 

As alluded to in Section~\ref{sec.SP} we cannot say in general whether $\Pi_{\theta,\theta'}(P)$ is the projection onto the closure of 
 $ E_{\theta,\theta'}(P)$ and along that of $E_{0}(P)\dotplus E_{\theta',\theta+2\pi}(P)$, but there are some important cases for which we can. 
 First, we have:  
 
 \begin{proposition}\label{prop:Appendix-CR.smoothingness}
 The following are equivalent:\smallskip 
   
 (i) For any $(x,\xi) \in T^{*}M\setminus 0$ there are no eigenvalues of $p_{m}(x,\xi)$ within $\overline{\Lambda_{\theta,\theta'}}$.\smallskip
 
 (ii) The sectorial projection $\Pi_{\theta,\theta'}(P)$ is a smoothing operator.\smallskip 
 
\noindent Moreover, if (i) and (ii) hold then $\Sp P \cap \Lambda_{\theta,\theta'}$ is finite and we have
    \begin{equation}
        \Pi_{\theta,\theta'}(P)=\sum_{\lambda \in \Sp P \cap \Lambda_{\theta,\theta'}} \Pi_{\lambda}(P).
    \end{equation}
Hence $\Pi_{\theta,\theta'}(P)$ has range $E_{\theta,\theta'}(P)$.
\end{proposition}
\begin{proof}
Since $\Pi_{\theta,\theta'}(P)$ is a (bounded) \psido\ projection, either it has order zero or it is smoothing. Thus $\Pi_{\theta,\theta'}(P)$ is a 
smoothing operator if, and only if, its zero'th order symbol is zero. By Proposition~\ref{prop:Appendix.sectorial-projection-PsiDO} 
the latter is the Riesz 
   projection $\Pi_{\theta,\theta'}(p_{m}(x,\xi))$ onto the root space associated to eigenvalues of $p_{m}(x,\xi)$ in $\Lambda_{\theta,\theta'}$. 
   Therefore $\Pi_{\theta,\theta'}(P)$ is smoothing if, and only if,  for any $(x,\xi) \in T^{*}M\setminus 0$ there are no eigenvalues of $p_{m}(x,\xi)$ within 
$\overline{\Lambda_{\theta,\theta'}}$. 

Assume now that for any $(x,\xi) \in T^{*}M\setminus 0$ there are no eigenvalues of $p_{m}(x,\xi)$ within 
$\overline{\Lambda_{\theta,\theta'}}$. Then there is an open angular sector $\Lambda$ containing $\overline{\Lambda_{\theta,\theta'}}\setminus 0$ such that 
  no eigenvalue of $p_{m}(x,\xi)$ lies in $\Lambda$. Then~(\ref{eq:NCR.minimal-growth}) tells us that $\Sp P\cap \Lambda_{\theta,\theta'}$ is finite and 
  for $R$ large enough there exists $C_{R\theta\theta'}>0$ such that we have
\begin{equation}
    \|(P-\lambda)^{-1}\|_{\cL(L^{2}(M,\cE))}\leq C_{R\theta\theta'}|\lambda|^{-1}, \qquad \lambda \in \overline{\Lambda_{\theta,\theta'}}, \quad 
    |\lambda|\geq R.
\end{equation}
It follows that in~(\ref{eq:SP.SAPi}) we may replace the integration contour $\Gamma_{\theta,\theta'}$ by a 
\emph{bounded} smooth contour $\Gamma$ which has index $-1$ and enlaces $\Sp P\cap \Lambda_{\theta,\theta'}$ but not the origin. Therefore, 
using~(\ref{eq:Appendix.proof.trick})  we see that $ \Pi_{\theta,\theta'}(P)$ is equal to
\begin{equation}
    \frac{1}{2i\pi}\int_{\Gamma} \frac{P}{\lambda(P-\lambda)} d\lambda = \!\!\! \sum_{\mu \in \Sp P\cap \Lambda_{\theta,\theta'}} \!\!\!
    \frac{-1}{2i\pi} \int_{\Gamma_{\mu}}\frac{P}{\lambda(P-\lambda)}d\lambda  = \!\!\!\sum_{\mu \in \Sp P\cap \Lambda_{\theta,\theta'}} \!\!\! \Pi_{\mu}(P).
\end{equation}
The proof is thus achieved. 
\end{proof}

Next, recall that $P$ is said to have \emph{a complete system of root vectors} when the total root space $\dotplus_{\lambda \in \Sp P} E_{\lambda}(P) $ is dense 
in $L^{2}(M,\cE)$. 

\begin{proposition}\label{prop:AppendixC.complete-system-root-vector}
    If $P$ has a complete system of root vectors then $\Pi_{\theta,\theta'}(P)$ is the projection onto $\overline{E_{\theta,\theta'}(P)}$   
    and along $E_{0}(P)\dotplus \overline{E_{\theta',\theta+2\pi}(P)}$.
\end{proposition}
\begin{proof}
    Let us first prove that $\ran \Pi_{\theta,\theta'}(P)$ is equal to $\overline{E_{\theta,\theta'}(P)}$. 
    We already know that the latter is contained in the former. Conversely, let $\xi$ be  in $\ran 
    \Pi_{\theta,\theta'}(P)$, so that $\Pi_{\theta,\theta'}(P)\xi=\xi$. Since $P$ has a complete system of root vectors there exists a sequence 
    $(\xi_{k})_{k\geq 0} \subset \dotplus_{\lambda \in \Sp P} E_{\lambda}(P) $ which converges to $\xi$
in $L^{2}(M,\cE)$. As $\xi_{k}$ is the sum of finitely many root vectors we have $\xi_{k}= \sum_{\lambda \in \Sp 
P}\Pi_{\lambda}(P)\xi_{k}$, where the sum is actually finite. Combining this with~(\ref{eq:AppendixC.Pitt'Pil}) gives
\begin{equation}
    \Pi_{\theta,\theta'}(P)\xi_{k}=\sum_{\lambda \in \Sp P}\Pi_{\theta,\theta'}(P)\Pi_{\lambda}(P)\xi_{k}=
    \sum_{\lambda \in \Sp P \cap \Lambda_{\theta,\theta'}}\Pi_{\lambda}(P)\xi_{k},
\end{equation}
so that $\Pi_{\theta,\theta'}\xi_{k}$ belongs to $E_{\theta,\theta'}(P)$. Since $\xi=\Pi_{\theta,\theta'}(P)=\lim_{k\rightarrow \infty}\Pi_{\theta,\theta'}\xi_{k}$ 
it follows that $\xi$ is in the closure of $E_{\theta,\theta'}(P)$. 
Hence $\ran \Pi_{\theta,\theta'}(P)=\overline{E_{\theta,\theta'}(P)}$. 
    
    Similarly, the projection $\Pi_{\theta',\theta+2\pi}(P)$ has range $\overline{E_{\theta',\theta+2\pi}(P)}$. Observe also that as 
    in~(\ref{eq:SAZF.decomposing-Wodzicki-projections})  we have 
    $\Pi_{\theta,\theta'}(P)+\Pi_{\theta',\theta+2\pi}(P)=1-\Pi_{0}(P)$. Thus,
    \begin{equation}
       \ran (1-\Pi_{\theta,\theta'}(P))= \ran \Pi_{0}(P) +\ran \Pi_{\theta',\theta+2\pi}(P)= E_{0}(P)\dotplus \overline{E_{\theta',\theta+2\pi}(P)}.
    \end{equation}
   Hence $\Pi_{\theta,\theta'}(P)$ is the projection onto $\overline{E_{\theta,\theta'}(P)}$ and along $E_{0}(P)\dotplus \overline{E_{\theta',\theta+2\pi}(P)}$.
\end{proof}

\begin{proposition}\label{prop:SP.normal}
    If $P$ is normal then $\Pi_{\theta,\theta'}(P)$ is the orthogonal projection onto 
    $\oplus_{\lambda \in \Sp P\cap \Lambda_{\theta,\theta'}} \ker (P-\lambda)$, where $\oplus$ denotes the Hilbertian direct sum on $L^{2}(M,\cE)$.
\end{proposition}
\begin{proof}
   Since $P$ is normal it diagonalizes on a Hilbert basis, that is, we have 
   \begin{equation}
       L^{2}(M,\cE)=\oplus_{\lambda\in \Sp P} \ker(P-\lambda),
        \label{eq:Appendix-spectral.orthogonal-splitting}
   \end{equation}
    where $\oplus$ denotes the Hilbertian direct sum on $L^{2}(M,\cE)$ (see~\cite[Thm.~V.2.10]{Ka:PTLO}). 
 In particular, $P$ has a complete system of root vectors, so by Proposition~\ref{prop:AppendixC.complete-system-root-vector}
 the sectorial projection $\Pi_{\theta,\theta'}(P)$ projects onto $\overline{E_{\theta,\theta'}(P)}$ and along $E_{0}(P)\dotplus \overline{E_{\theta',\theta+2\pi}(P)}$. 
 
 On the other hand, the orthogonal decomposition~(\ref{eq:Appendix-spectral.orthogonal-splitting}) 
 implies that for every $\lambda\in \Sp P$ we have $E_{\lambda}(P)=\ker (P-\lambda)=\ker (P^{*}-\overline{\lambda})$. Thus,
\begin{equation}
    \overline{E_{\theta,\theta'}(P)}=\oplus_{\lambda \in \Sp P\cap \Lambda_{\theta,\theta'}} \ker(P-\lambda).
\end{equation}
Similarly $E_{0}(P)\dotplus \overline{E_{\theta',\theta+2\pi}(P)}$ is equal to 
\begin{equation}
    \ker P\oplus[\oplus_{\lambda \in \Sp P\cap \Lambda_{\theta',\theta+2\pi}} \ker 
      (P-\lambda)]=E_{\theta,\theta'}(P)^{\perp}.
\end{equation}
Hence $\Pi_{\theta,\theta'}(P)$ is the orthogonal projection onto  $\oplus_{\lambda \in \Sp P\cap \Lambda_{\theta,\theta'}} \ker (P-\lambda)$.
\end{proof}
As an immediate consequence we get: 

\begin{corollary}
    When $P$ is selfadjoint the sectorial projection $\Pi_{\uparrow\downarrow}(P)$ is the orthogonal projection onto the negative eigenspace of $P$. 
\end{corollary}

There are well known examples due to Seeley~\cite{Se:SESPDO} and Agranovich-Markus~\cite{AM:SPEPDOFSO} 
of elliptic differential operators without 
a complete system root vectors. In these examples the  principal 
symbol does not admit a spectral cutting. However, even when the principal symbol does admit a spectral cutting the best positive result about completeness 
result seems to be the following.

\begin{proposition}[{\cite[Thm.~3.2]{Ag:EEGEBVP}}, {\cite[Appendix]{Bu:FPEDO}}, 
{\cite[Thm.~6.4.3]{Ag:EOCM}}]\label{prop:SP.complete-root-vector-system}
    Assume that the principal symbol of $P$ admits spectral cuttings $L_{\theta_{1}},\ldots, L_{\theta_{N}}$  dividing the complex planes into 
angular sectors of apertures~$<\frac{2n\pi}{m}$. Then the system of root vectors of $P$ is complete. 
    \end{proposition}

This result follows from a criterion due to Dunford-Schwartz~\cite[Cor.~XI.9.31]{DS:LOII} for closed operators on a Hilbert spaces with a resolvent in some 
Schatten ideal. Combining it with Proposition~\ref{prop:AppendixC.complete-system-root-vector} thus gives: 

\begin{proposition}\label{prop:SP.DS-criterion} If the principal symbol of $P$ admits spectral cuttings dividing the complex plane into 
angular sectors of apertures~$<\frac{2n\pi}{m}$, then $\Pi_{\theta,\theta'}(P)$ is the projection onto $\overline{E_{\theta,\theta'}(P)}$ and along 
$E_{0}(P)\dotplus \overline{E_{\theta',\theta+2\pi}(P)}$. 
\end{proposition}

Finally, if we only want to determine the range of $\Pi_{\theta,\theta'}(P)$ then we have:
\begin{proposition}\label{prop:SP.range-criterion}
  If the principal symbol of $P$ admits spectral cuttings $L_{\theta_{1}},\ldots, L_{\theta_{N}}$ dividing the angular sector $\Lambda_{\theta,\theta'}$ into 
angular sectors of apertures~$<\frac{2n\pi}{m}$. Then the range of $\Pi_{\theta,\theta'}(P)$ is equal to $\overline{E_{\theta,\theta'}(P)}$. 
\end{proposition}
\begin{proof}
    The operator $\tilde{P}$ induced by $P$ on $\ran \Pi_{\theta,\theta'}(P)$ has spectrum $\Sp P\cap \Lambda_{\theta,\theta'}$ and its resolvent is 
    also in the Schatten ideal $\cL^{\frac{m}{n}+\epsilon}$ for any  $\epsilon >0$. 
   Moreover,  the condition on the principal symbol implies that $P$ has finitely many rays of minimal growth $L_{\theta_{1}'},\ldots, L_{\theta_{N}'}$ 
    dividing $\Lambda_{\theta,\theta'}$ into angular 
sectors of aperture~$<\frac{2n\pi}{m}$. Henceforth
$\tilde{P}$ admits a finite sequence of rays of minimal growth dividing $\C$  into angular 
sectors of aperture~$<\frac{2n\pi}{m}$. It then follows from~\cite[Cor.~XI.9.31]{DS:LOII} that the total root space of $\tilde{P}$,  that is, $E_{\theta,\theta'}(P)$, is 
dense in $\ran \Pi_{\theta,\theta'}(P)$. 
\end{proof}
  \begin{acknowledgements} 
I am indebted to 
Maxim Braverman, Sasha Gorokhovsky, Gerd Grubb, Boris Mityagin, Henri Moscovici, Victor Nistor and Mariusz Wodzicki for helpful and 
 stimulating discussions and to Sasha Dynin for having provided me with a copy of~\cite{Wo:CWSP}. 
 In addition, I also wish to thank for their hospitality the Max Planck Institute for Mathematics (Bonn, Germany)
 and the IH\'ES (Bures-sur-Yvette, France) where parts of this paper were written.
\end{acknowledgements}

\end{document}